\documentclass[11pt]{amsart}
\usepackage{xcolor}
\usepackage{graphicx}
\usepackage{latexsym}
\usepackage{amsfonts,amsmath,amssymb}
\newtheorem{theorem}{Theorem}[section]
\newtheorem{proposition}[theorem]{Proposition}
\newtheorem{lemma}[theorem]{Lemma}

\newtheorem{corollary}[theorem]{Corollary}
\newtheorem{example}[theorem]{Example}

\usepackage{color}
\def\si{\par\smallskip\noindent}
\def\bi{\par\bigskip\noindent}

\def\de{\delta}
\def\eps{\varepsilon} 
\def\la{\lambda}
\def\a{\alpha}
\def\be{\beta}
\def\de{\Delta}

\def\ga{\gamma}
\def\part{\partial}

\def\bk{\bold k}
\def\bv{\bold v}
\def\b1{\bold 1}

\newcommand{\beq}{\begin{equation}}
\newcommand{\eeq}{\end{equation}}

\theoremstyle{remark}

\numberwithin{equation}{section}
\date{\today}

\begin{document}

\title[Plane increasing trees]{Random plane increasing trees: asymptotic enumeration of vertices by distance from leaves}

\author{Mikl\'os B\'ona}
\address{Department of Mathematics, University of Florida, $358$ Little Hall, PO Box $118105$,
Gainesville, FL, $32611-8105$ (USA)}
\email{bona@ufl.edu}
\author{Boris Pittel}
\address{Department of Mathematics, The Ohio State University, $231$ West $18$-th Avenue, Columbus, Ohio $43210-1175$ (USA)}
\email{bgp@math.ohio-state.edu}

\keywords
{search tree,  root, leaves, ranks, enumeration, asymptotic, distribution,
numerical data}

\subjclass[2010] {05A05, 05A15, 05A16, 05C05, 06B05, 05C80, 05D40, 60C05}

\begin{abstract}
We prove that for any fixed $k$, the probability that a random vertex of a random increasing plane tree 
is of rank $k$,  that is, the probability that a random vertex is at distance $k$ from the leaves, converges to a constant $c_k$ as the size $n$ of the tree goes to infinity. 
{\color{blue} We prove that  $1-\sum_{j\le k} c_k<\tfrac{3^{k+1}}{(2k+1)!}$, so that the tail of the limiting rank
distribution is super-exponentially narrow. We prove that the latter property holds uniformly for all finite $n$ as well.}  More generally,
we prove that %a similar
%result for randomly selected $r$-tuples of vertices that shows that 
the ranks of a finite uniformly random set
of vertices are asymptotically independent, each with distribution $\{c_k\}$. We compute the exact value of $c_k$ for $0\leq k\leq 3$, demonstrating that
the limiting expected fraction of vertices with rank $\le 3$ is $0.9997\dots$. We show that with probability $1-n^{-0.99\eps}$  the highest rank of a vertex in the tree is sandwiched between 
$(1-\eps)\log n /\log\log n$ and $(1.5+\eps)\log n/\log\log n$, {\color{blue} and that this rank is asymptotic to $\log n/\log\log n$ with probability $1-o(1)$.}
\end{abstract}
 
\dedicatory{We dedicate this paper to the precious memories of Mikl\'os B\'ona, the  father of the first author,
 and of Irina Pittel, the wife of the second author.}

\maketitle

\section{Introduction}
\subsection{Definitions and Background}
The analysis of various parameters of random rooted trees is a classic subject if those parameters relate to a near-root part of a tree or to 
global tree structure. See \cite{flajolet} for a comprehensive treatment of these results. A more recent area of interest is the random tree
fringe, that is, the part close to the leaves. A list of articles on this subject can be found in \cite{bona-pittel}. 

A particularly interesting problem is the enumeration of vertices of a given {\em rank}. 
If $v$ is a vertex of a  tree $T$, then let the {\em rank} of $v$ be the number of edges in the
shortest  path from $v$ to a leaf of $T$ that is a {\em descendant} of $v$. So leaves are of rank 0, parents of leaves are of rank 1, and so on.
Let us fix a variety of trees, such as binary search trees, or non-plane binary trees, or increasing trees. Let $n$ be a fixed positive integer, 
and let us compute the probability that a randomly selected vertex of a randomly selected tree of size $n$ is of rank $k$. What can we say about this probability as $n$ goes to infinity? 

This question for various models of random trees has been the subject of vigorous research in the last few years. 
Perhaps the first was a random binary search tree process $\{T_n\}$: conditioned on a binary tree $T_n$ with $n$ vertices, $T_{n+1}$ is obtained by selecting uniformly at random (u.a.r.) one of $n$ available positions for a new leaf and putting the next vertex into it. It was proved in  \cite{protected} that for any fixed $k$, the probability in question converges to a rational number $c_k$, and $c_k$ was computed for $k\leq 3$. Then in \cite{bona-pittel}, the present authors computed $c_4$ and $c_5$, and proved lower and upper bounds for 
$c_k$ in general. They proved that as $k$ grows, $c_k$ exponentially decays. They also showed that  the ranks of the uniformly random, fixed size, sample of vertices are asymptotically independent. For increasing non-plane trees where the
next vertex is attached to a previous vertex selected u.a.r. from the existing vertices, known also as recursive trees, ranks were extensively studied by Holmgren and Janson
in \cite{hj}. They also studied ternary trees from this aspect \cite{Holmgren}. In \cite{bona-mezo}, the first author and Istv\'an Mez\H{o} proved the existence of the limits $c_k$ for non-plane 1-2 trees and plane 1-2-trees. 

A natural tree variety that is missing from the above list is {\em plane increasing trees}. These are rooted plane trees  whose vertices are bijectively labeled with the elements of $[n]$ so that each vertex has a 
{\em larger} label than its parent, and the children of every parent form an {\bf ordered} set. This tree variety is the subject of our paper. As these are the only kind of trees that will appear in this 
work, we will often simply call these rooted  plane increasing trees just {\em trees}. These trees have been studied under various names, such as {\em plane-oriented trees} \cite{Hwang,Mahmoud}, {\em PORT}s (see page 14 of \cite{book-Drmota}), 
or scale free trees \cite{Bar}. 

Given a tree with $n-1$ vertices of {\it total\/} degrees $d_1,\dots, d_{n-1}$, there are
$d_{\text{root}}+1+ \sum_{j\in [2,n-1]} d_j =2n-3$ different ways to attach vertex $n$ to a vertex in $[n-1]$. {\color{blue} Indeed, there are $d_{\text{root}}+1$ ways to attach $n$ to the root, and $d_j$ ways to attach $n$ to any other non-root vertex $j\in [2,n-1]$. So the total number of these trees
is $(2n-3)!!:=1\cdot 3 \cdot \cdots \cdot (2n-3)$.} Our tree $T_n$ is chosen u.a.r. from the set of all $(2n-3)!!$ trees. We can view $T_1, T_2,\dots$ as snapshots of the tree
growing (Markov) process: conditioned on $T_{n-1}$, vertex $n$ becomes a child of a vertex $j\in [n-1]$ taking up one of 
$d_j$ 
possible positions in the ordered set of the {\color{blue} other children of $j$}. So, unlike in recursive trees, a host vertex for the vertex $n$ is chosen from existing vertices with probability proportional to total degree of {\color{blue} that} potential host. The second author \cite{Pit1} proved that {\color{blue} the
height, that is, the} length of the longest path from the root to a leaf of this tree, scaled by $\log n$ converges, almost surely, to $1.795560738$, which is the root of the equation $2y/e=\exp(1/2y)$. Later, Drmota \cite{Drmota} extended this result to a broad variety of increasing
trees introduced by Bergeron, Flajolet, and Salvy \cite{Berg}.

In a different direction, $\{T_n\}_{n\ge 1}$ is a well-known case of the random tree (graph, if cycles allowed) process with ``preferential attachment''. See Barab\'asi and Albert \cite {Bar}, M\'ori \cite{Mor1}, \cite{Mor2}, Bollob\'as and Riordan \cite{BolRio}, 
Bollob\'as \cite{Bollobas}, van der Hofstad \cite{Hof}, Frieze and Karo\'nski \cite{Frieze},  Pittel \cite{Pit2}.

\subsection{Our Main Results} \label{mainresults}
As we mentioned above, in this paper, the tree variety we study is that of increasing plane trees. In Section \ref{asymptotics}, we prove that for this tree variety, the numbers $c_k$ exist for all $k$, and that $c_k$ is equal to the value of a certain integral. We provide the explicit value of this integral for $k=0$ and $k=1$, and its approximate value (to eleven decimals) for 
$k=2$ and $k=3$. {\color{blue} {\color{blue} We also prove a non-obvious result} that $\sum_{k\geq 0} c_k =1$, so the numbers $c_k$ form a probability 
distribution.

 In order to see that this fact is not obvious and has to be proved, consider the trivial variety of trees in which every non-leaf vertex has exactly one child. 
In that case, the ratio of vertices of rank $k$ in a tree of size $n$ is $c_k(n)=1/n$ if $k<n$, so $c_k=0$ for any fixed $k$. That is,  for this tree variety, $\sum_k c_k =0 <1$. 

For an arbitrary tree variety, only the following holds.  As $\sum_{k=0}^n c_k(n) =1$, it follows that
 for each fixed bound $B$,  we have
$\sum_{k\le B} c_k(n)\le 1$,
and taking limits as $n$ goes to infinity, 
$\sum_{k=0}^B c_k \le 1$. Now taking limits as $B$ goes to infinity, 
we get
$\sum_{k\ge 0} c_k \le 1$. So the sum of the $c_k$ is never more than 1, i.e. $\sum_{k\ge 0}c_k\in [0,1]$. To show for a particular increasing tree that  $\sum_{k\geq 0} c_k=1$, one has to prove tightness condition, namely that $\sum_{k=K}^n c_k(n)$ can be made arbitrarily small for all $n$ by selecting a {\em constant} $K$ that is sufficiently large. Less formally,  one has to show that, as a tree in question grows indefinitely, with high probability no
substantial portion of the vertices will have higher and higher ranks, exceeding any fixed $k$. That's what we will prove 
for the tree growing process studied in this paper. }

In fact, we will prove not only that $\sum_k c_k=1$, but also the stronger statement that
 $1-\sum_{j\le k}c_j\le \tfrac{3^{k+1}}{(2k+1)!}$. {\color{blue} So the tail of the limiting rank distribution is super-exponentially narrow. 
The latter property holds uniformly for all finite $n$ as well.} These results together imply that as $n$ goes to infinity,  the limiting probability that a random vertex of a random tree 
has rank four or more is approximately 0.0002843360. So about 99.97 percent of all vertices have rank three or less. 

Using our results on the expected number of vertices of a given rank, we are able to prove lower and upper bounds on the {\em highest rank} $R_n$ of a vertex in a tree of size $n$ selected 
uniformly at random. {\color{blue} First, in Theorem \ref{thm1}, we show that, for each $\eps>0$, the inequality $R_n \leq  (1.5+\epsilon) \log n / \log\log n$ holds with probability $1-n^{-\eps-o(1)}$. (Bollob\'as and Riordan \cite{Bollobas} proved that with probability $1-o(1)$, the diameter of the preferential attachment tree is asymptotic to $\log n/\log\log n$).
Later, in Theorem \ref{thm4}, focusing on vertices whose descendant tree is a path we prove that if $0<\epsilon< 1$, then with probability $1-n^{-\epsilon}$, there exist $\Theta (n^\epsilon)$ vertices with
rank $k=(1-\epsilon)\log n/\log \log n$. (This result obviously implies Bollob\'as-Riordan's lower bound for the diameter.) In other words, with probability $1-o(1)$, the highest rank is asymptotic to $\log n/ \log \log n$.}
%and it is asymptotic to $\log n/\log\log n$ with probability $1-o(1)$.}

In Section \ref{rtuples}, we turn to the counting of $r$-tuples of vertices with given ranks. We are able to prove that the ranks of a finite, uniformly random, set of vertices
are asymptotically independent. In other words, if {\color{blue} $R_n(1),R_n(2),\cdots ,R_n(r)$} denote the ranks of $r$ vertices, chosen uniformly at random, with order and without replacement, from the
set of $n$ vertices of a random tree, then 

\[
\Bbb  P\left ( \bigcap_{i=1}^r ( R_n(i)=k_i)\right)=(1+O(n^{-1})) \prod_{j=1}^r c_{k_j}.
\]

\section{Counting vertices by ranks via generating functions}\label{Sect1}

\subsection{Individual vertices and tuples of vertices}\label{Sub1.1}
{\color{blue} The following result is well-known \cite{Mahmoud}, and we have proved it in a purely combinatorial way in the introduction. Nevertheless, here we include another proof,  to show the method we will use later in more complicated circumstances. }

\begin{proposition} \label{countingtrees}
Let $N\geq 2$. Then the number $t_n$ of all rooted increasing plane trees on vertex set $[n]$ is given
by the formula $t(n)=(2n-3)!!$.
\end{proposition}

\begin{proof} Note that  $t(1)=1$ and for $n\ge 2$ we have the recurrence 
\begin{equation}\label{new-9}
t(n)=\sum_{s\ge 1}\sum_{j_1+\cdots+j_s=n-1}\frac{(n-1)!}{j_1!\cdots j_s!}\prod_{i=1}^s t(j_s),\quad (t(0):=0),
\end{equation}
where $s$ is the number of the children of the root.   Multiplying both sides by $z^{n-1}/(n-1)!$, summing over $n\ge 2$ and adding $z\,t(1)$, we obtain that the power series $T(z):=\sum_{n\ge 1}\frac{t(n)}{n!} z^n$  satisfies the equation 
\begin{equation} \label{treeeq}
T'(z)=\sum_{j\ge 0} T^j(z)= (1-T(z))^{-1},\quad T(0)=0,
\end{equation}
if $|T(z)|<1$. The solution of this differential equation is $T(z)= 1-\sqrt{1-2z}$, where we use the main branch of the square root,
and $|T(z)|<1$ if $|z|<1/2$. It follows that 
\[
\frac{t(n)}{n!}=-[z^n](1-2z)^{1/2}=\frac{(2n-3)!!}{n!}\longrightarrow t(n)=(2n-3)!!.
\]
%i.e. $t(n)=(2n-3)!!$, and $|T(z)|\le T(|z|) <1$ if $|z|<1/2$.
\end{proof}

{\bf Note.\/} {\color{blue} We could have used the symbolic method (Flajolet and Sedgewick \cite{flajolet}) and obtain (\ref{treeeq}) directly, but we decided to present an explicit, admittedly pedestrian, argument where we start with tree-based recurrence for the tree counts in question and then derive recurrences for
generating functions. This short proof will serve as a self-contained template for our considerably more technical proofs coming next, especially when we will turn to counting tuples of vertices of given ranks.} 

{\bf Notation.\/} Suppose $|z|<1/2$. Let $A_k(z)=\sum_{n\geq 1}a_k(n) z^n/n!$ be the (exponential) generating function for the number of  all vertices of rank $k$ in all rooted plane trees on vertex set $[n]$.
Let $B_k(z)=\sum_{n\geq 1}b_k(n) z^n/n!$ be the generating function for the number of  trees in which {\em the root} is of rank $k$, so  we have $B_0(z)=z$. 
%Clearly $a_k(n),\,b_k(n)\le t(n)$, which we express as $A_k(z),\,B_k(z)\le T(z)$. 
We will also need
$B_{\geq k}(z)$, the generating function  for the number of trees in which the root is of rank $\ge k$. So we have$B_{\ge 0}(z)=T(z)$, and
$B_{\ge 1}(z)$ is the generating function of numbers of trees with at least two vertices, i.e. $B_{\ge 1}(z)=T(z)-z$. 
We will use the {\color{blue} inequality}
$B_{\ge k+1}(z)\le B_{\ge k}(z)$ to express the fact  that for all $n$, the coefficient  $[z^n]B_{\ge k}(z)$ decreases as $k$ increases.

\begin{proposition} \label{recur} For all $k$, the generating functions \label{general} $A_k(z)$ and $B_k(z)$ satisfy
\begin{equation}\label{r=1}
 \begin{aligned}
 A_k'(z)=\frac{A_k(z)}{(1-T(z))^2}& +B_k'(z),\,\,B_k(z)=B_{\ge k}(z)-B_{\ge k+1}(z),\,\, k\ge 0,\\
 B_{\geq k}'(z) =& \frac{B_{\ge k-1}(z)}{1-B_{\geq k-1}(z)}, \,\, k\ge 1,\,\,B_{\ge 0}(z)=T(z).
\end{aligned}
\end{equation}
\end{proposition}

{\color{blue}
\begin{proof} Let $X_k(n)$ stand for the total number of rank $k$ vertices in the uniformly random tree, so that $\Bbb E[X_k(n)]=\tfrac{a_k(n)}{t(n)}$. Let $p_k(n)$, (resp. $p_{\ge k}(n)$) stand for probability that the root $1$ is of rank $k$ (resp. rank $\ge k$). {\color{blue} So, the formulas $b_k(n)=t(n)p_k(n)$ and $b_{\ge k}(n)=t(n)p_{\ge k}(n)$ provide the total number of trees on $[n]$ with root of rank $k$ and of 
rank $\ge k$,  respectively.} 
For $k\ge 1$, the root is of rank at least $k$ if and only if the root has children and each of them is of rank at least $k-1$. Therefore, for $n\ge 2$,
\begin{equation}\label{n1}
p_{\ge k}(n)=\tfrac{(n-1)!}{t(n)}\sum_{s\ge 1}\sum_{j_1+\dots+j_s=n-1}\prod_{r=1}^s\tfrac{t(j_r)\,p_{\ge k-1}(j_r)}{j_s!},\quad j_1,\dots,j_s\ge 1.
\end{equation}
So, $B_{\ge k}(z)$ and $B_k(z)$, the exponential generating functions of $\{b_{\ge k}(n)\}$ and $\{b_k(n)\}$, satisfy
\begin{align*}
B^\prime_{\ge k}(z)&=\sum_{s\ge 1} (B_{\ge k-1}(z))^s=\tfrac{B_{\ge k-1}(z)}{1-B_{\ge k-1}(z)}, \,\,k\ge 1;\quad B_{\ge 0}(z)=T(z),\\
B_k(z)&=B_{\ge k}(z)-B_{\ge k+1}(z),\quad k\ge 0.
\end{align*}
{\color{blue} The  series  $B^\prime_{\ge k}(z)$ above converges for $|z|<1/2$ and $k\ge 1$ because
\[
|B_{\ge k-1}(z))|\le B_{\ge k-1}(|z|)\le B_{\ge 0}(|z|)=T(|z|)<T(1/2)=1.
\]

Let us turn to $A_k(z)$. Now, in the random tree $X_k(n)$ equals $\Bbb I(\text{root is of rank }k)$ plus the total sum of all vertices of rank $k$ in the subtrees rooted at the children of the root. } The probability that the root has $s$ children with  ordered subtrees of sizes $j_1,\dots,j_s$ is
$
\tfrac{(n-1)!}{t(n)}\prod_{r=1}^s\tfrac{t(j_r)}{j_r!}.
$
So, by linearity of expectation, for $n\ge 2$, the following holds. 
\begin{multline}\label{??}
\Bbb E[X_k(n)]=p_k(n)\\
+\tfrac{(n-1)!}{t(n)}\sum_{s\ge 1}\sum_{j_1+\dots+j_s=n-1}\prod_{r=1}^s\tfrac{t(j_r)}{j_r!}\biggl(\sum_{r'=1}^s
\Bbb E[X_{k}(j_{r'})]\biggr).
\end{multline}
Multiplying both sides by $\tfrac{t(n)z^{n-1}}{(n-1)!}$, and summing over $n\ge 2$, we obtain
\begin{multline*}
\tfrac{d}{dz}\sum_{n\ge 1}\tfrac{t(n)\Bbb E[X_k(n)]}{n!}\,z^n=\tfrac{d}{dz}\sum_{n\ge 1}\tfrac{t(n)p_k(n)}{n!}\,z^n\\
+\sum_{s\ge 1}\sum_{r=1}^s\sum_{j_{r}\ge 1}\tfrac{t(j_r)E[X_k(j_r)]}{j_r!}\,z^{j_r}\sum_{j_{r'}\ge 1,\,r' \neq r}\,\,\prod_{r"\neq r}\tfrac{t(j_{r")}}{j_{r"}!}\,z^{j_{r"}},
\end{multline*}
or equivalently
\[
A'_k(z)=B'_k(z)+A_k(z)\sum_{s\ge 1}s T^{s-1}(z)=B'_k(z)+\tfrac{A_k(z)}{(1-T(z))^2}.
\]
\end{proof} }

\begin{example}\label{Ex1} {\color{blue} Let $n\geq 2$. Then the number of all leaves in all trees of size $n$ is $(2n-1)!!/3$. }   Indeed, by
Proposition \ref{general} with $k=0$ and $B_0(z)=z$, we have 
\[
A_0(z)= \frac{2z-1}{3} + \frac{1}{3\sqrt{1-2z}}\Longrightarrow \frac{a_0(n)}{n!}=\frac{1}{3}\frac{(2n-1)!!}{n!}.
\]
\end{example}

{\color{blue}
\subsubsection{Counting ordered tuples of vertices by ranks}  
%An {\em ordered partition} of the $n$-element set $[n]=\{1,2,\cdots, n\}$ is a sequence of disjoint  non-empty sets called {\em blocks} whose union
%is  $[n]$. Note that the order of the various blocks matters, but the order of the elements within the blocks does not. For instance, $\{(1,3),(2)\}$ and $\{(3,1),(2)\}$ are identical as
%ordered partitions, but $\{(2),(1,3)\}$ is different. As the order of entries within each block does not matter, we are writing them in increasing order by convention. 

%In what follows, we will need the concept of an ordered partition of a generic $r$-tuple of nonnegative integers, instead of the set $[n]$. Otherwise, the concept is analogous to what
%we described in the preceding paragraph. An ordered partition of an $r$-tuple $R$ of nonnegative integers is a sequence of non-empty blocks so that each block is a sub-tuple of $R$, and
%the {\em disjoint} union of the subtuples is $R$. For instance, if $R=(2,3,1,3,4)$, then $\{(2,1,3),(3,4)\}$ is an ordered partition of $R$, where the block $(2,1,3)$ contains the first, third, and 
%fourth entries of $R$, and the block $(3,4)$ contains the second and fifth entries of $R$. As the entries of $R$ are not necessarily disjoint, it may be helpful to think about the occurrences of entries appearing with
%multiplicity higher than one as appearances of different color. }
Given an ordered $r$-tuple $\bk=(k_1,\dots,k_r)$ of (non-negative) integers, let $a_{\bk}(n)$ denote
the total number of occurrences of ordered $r$-tuples of distinct vertices $\bv=(v_1,\dots,v_r)$ with ranks $\bk$ in all  trees on $[n]$. (Proposition \ref{recur} {\color{blue} deals with a} special case $r=1$.)
Given $k_1$ and an $(r-1)$-long  ordered tuple $\bk'=(k_2,\dots,k_r)$, let $b_{\bk}(n)$, ($\bk:=(k_1,\bk')$), denote the total number of occurrences of $r$-long  tuples $(root, \bv')$ such that the root has rank $k_1$ and $r-1$ distinct vertices $\bv'$ have ranks $\bk'$. Let $A_{\bk}(z)$ and
$B_{\bk}(z)$ denote the respective generating functions. An {\em ordered partition} of a set or $r$-tuple is a partition of that
set or $r$-tuple into nonempty blocks in which {\em the set of blocks is linearly ordered}. There is no new ordering on each 
individual block; the individual blocks inherit the original ordering.

\begin{proposition}\label{recur,gen} Let us assume that  $r>1$ and $\bk$ is an ordered $r$-tuple. Let $(\bk^{(1)},\dots,\bk^{(\tau)})$ denote an ordered partition of $\bk$ into $\tau\le r$ non-empty sub-tuples, each sub-tuple inheriting its ordering from $\bk$. Then 
\begin{align*}
A'_{\bk}(z)&= \sum_{\tau\le r}(1-T(z))^{-\tau-1}\biggl(\,\sum_{\bk^{(1)}\!,\dots,\bk^{(\tau)}\neq \emptyset}\,\,\prod_{u=1}^{\tau} A_{\bk^{(u)}}(z)\biggr)\\
%{|\bk^{(u)}|!}
&+\sum_{u\in [r]} B'_{k_u, (k_t)_{t\in [r]\setminus\{u\}}}(z).
\end{align*}
Further, with $\bk=(k_1,\bk')$, let $(\bk^{(1)},\dots, \bk^{(\tau)})$ denote a generic ordered partition of $\bk'$ into $\tau$ non-empty ordered sub-tuples. If $\bk'\neq\emptyset$, then
\begin{align*}
B'_{\bk}(z)&\le \sum_{\tau<r}\,\,(1-T(z))^{-\tau-1}\biggl(\,\sum_{\bk^{(1)}\!,\dots,\bk^{(\tau)}\neq \emptyset}\,\,\prod_{u=1}^{\tau} A_{\bk^{(u)}}(z) \biggr).
\end{align*}
This inequality means that, for each $\nu\ge 0$, $[z^{\nu}] B'_{\bk}(z)$ is at most the coefficient {\color{blue} of} $z^{\nu}$ in the generating function on the RHS.
%we have $B_{\bk}(z)\le B^+_{\ge k_1,\bk'}(z)$, where the $B^+_{\ge k_1,\bk'}(z)$ are power series defined 
%by the following recurrence.
%For $k_1\ge 0$ and $|\bk'|>0$, 
%\begin{multline*}%\label{new-17}
%\begin{aligned} 
%\frac{d B^+_{\ge k_1,\bk'}}{dz}\!\!=\!\!\sum_{s\ge 1}\frac{1}{(1-B_{\ge k_1-1}(z))^{s+1}}
%\times\sum_{\cup_{u=1}^{s}\bk^{(u)}=\bk';\atop \bk^{(u)}\neq \emptyset}\prod_{u=1}^s 
%B^+_{\ge k_1-1,\bk^{(u)}}(z).
%\end{multline*}
%\end{equation}
%For $|\bk'|=0$ we are back to $r=1$, and for $k_1\ge 0$, we have $B^+_{\ge k_1, \bk'}(z)=B_{\ge k_1}(z)$, defined at the bottom equation in Proposition \ref{general}.
\end{proposition}

\begin{example} 
Let $\bk=(k_1,k_2,k_3)=(3,7,7)$. Then $a_{\bk}(n)$ denotes
the total number of occurrences of ordered triples $\bv = \{v_1,v_2,v_3\}$ of vertices in all trees on $[n]$ so that $v_1$ is of rank 3, and $v_2$ and $v_3$ are both of rank $7$. 
Furthermore, $b_{\bk}(n)$ denotes the number of occurrences of all ordered triples $(root,v_2,v_3)$  in all such trees such that the root has rank 3, and $v_2,$ and $v_3$ are both of rank $7$. 

Then $\bk$ has $13$ ordered partitions into non-empty  sub-tuples, each inheriting its ordering from $\bk$:
 $(\{k_1,k_2,k_3\})$, $(\{k_1,k_2\},\{k_3\})$, 
$(\{k_3\}, \{k_1, k_2\})$, $(\{k_1,k_3\}, \{k_2\})$, $(\{k_2\}, \{k_1,k_3\})$, $(\{k_1\}, \{k_2,k_3\})$, 
$(\{k_2,k_3\}, \{k_1\})$, $(\{k_1\}, \{k_2\},$\linebreak $\{k_3\})$, $(\{k_2\}), \{k_1\}, \{k_3\})$, $(\{k_2\}, \{k_3\}, \{k_1\})$,
 $(\{k_3\}, \{k_1\}, \{k_2\})$,  
$(\{k_3\}, \{k_2\},$\linebreak $\{k_1\})$. Therefore, the first equality in Proposition \ref{recur,gen}
reduces to 
\begin{eqnarray*}
A_{(3,7,7)}'(z) & = & \tfrac{A_{(3,7,7)\}}(z)}{(1-T(z))^2} + \tfrac{2A_{(\{3,7\}, \{7\})}(z)+2A_{(\{7\},\{3,7)\})}(z)
+A_{\{(3), (7,7)\}}(z)}{(1-T(z))^3} \\
 & + &  B_{(\{3\}, \{7,7)\} }'(z)+2B'_{(\{7\}, \{3,7\})}(z).
\end{eqnarray*}
Further, $B_{\geq 3,(7,7)}(z)$ is the generating function for the numbers of ordered triples  $(root,v_2,v_3)$ so that the root has rank {\em at least three}, and $v_2$, $v_3$ are both of rank 7. So, the inequality in  Proposition \ref{recur,gen} reduces to 
\[
B_{(3,7,7)}'(z) \leq \tfrac{A_{(7,7)}(z)}{(1-T(z))^2} + \tfrac{A_7^2(z)}{(1-T(z))^3}.
\]
Likewise 
\[
B_{(\{7), \{3,7\}) }'(z)\le \tfrac{A_{(3,7)}(z)}{(1-T(z))^2}+\tfrac{2A_3(z) A_7(z)}{(1-T(z))^3}.
\]
\end{example}

\begin{proof} Let $X_{\bk}(n)$ stand for the total number of $r$-long tuples of distinct vertices of ranks $\bk=(k_1,\dots,k_r)$  in the random tree, so that $\Bbb E[X_{\bk}(n)]=\tfrac{a_{\bk}(n)}{t(n)}$.
Let $Y_{\bk}(n)$ (resp., $Y_{\ge k_1, \bk'}$) stand for the number of $r$-long tuples of distinct vertices $(1, v_2,\dots,v_r)=:(1,\bv')$ such that the root $1$ is of rank $k_1$ (of rank $\ge k_1$ resp.)  
and vertices forming $\bv'$ are of ranks $\bk':=(k_2,\dots,k_r)$. So, $\Bbb E[Y_{\bk}(n)]=\tfrac{b_{\bk}(n)}{t(n)}$ and $\Bbb E[Y_{\ge k_1, \bk'}(n)]=
\tfrac{b_{\ge k_1,\bk'}(n)}{t(n)}$. %.:=t(n)\Bbb E[Y_{\ge k_1, \bk'}(n)]$ are the total number of occurrences of tuples $\{1,\bv'\}$ in all trees on $[n]$, with root of rank $k_1$ and of rank $\ge k_1$ respectively, and $\bv'$ of ranks $\bk'$. For $k_1\ge 1$, the root is of rank at least $k_1$ if and only if the root has children and each of them is of rank at least $k_1-1$. 

Let us write a recurrence for $\Bbb E[X_{\bk}(n)]$. The $r$-long tuples $\bv$ that contain the root $1$ contribute  $\sum_{u\in [r]} \Bbb E\bigl[Y_{k_u, (k_t)_{t\in [r]\setminus\{u\}}}\bigr]$ to this expectation. (The $u$-th term corresponds to the case when $k_u$, i.e. the $u$-th component of $\bk$, is the root's rank.) Consider the tuples $\bv$ that do not contain the root $1$. Let $s\ge 1$ be a generic number of the root's ordered children, and let $[r]$ be a disjoint, union of $s$ subsets, some of which may be empty. This ordered partition of $[r]$ determines the  partition of the ordered tuple $\bk$ into $s$
sub-tuples $\bk^{(1)},\dots, \bk^{(s)}$ of total length $r$, where each sub-tuple inheriting the order of ${\bf k}$. Conditioned on $s$ and the sizes $j_1,\dots,j_s$ of the subtrees rooted at $s$ children of the root $1$, the subtrees are independent. Therefore, the conditional expected number of $r$-long tuples of rank $\bk$, composed of $s$ sub-tuples of ranks $\bk^{(1)},\dots, \bk^{(s)}$,
is $\prod_{u=1}^s\Bbb E[X_{\bk^{(u)}}(j_u)]$.  Here $\Bbb E[X_{\bk^{(u)}}(j_u)]:=1$ if $\bk^{(u)}=\emptyset$. Therefore, for $n\ge 2$, we have
\begin{eqnarray*}
\Bbb E[X_{\bk}(n)] & =  & \sum_{u\in [r]} \Bbb E\bigl[Y_{k_u, (k_t)_{t\in [r]\setminus\{u\}}}(n)\bigr]\\
&+&  \tfrac{(n-1)!}{t(n)}\sum_{s\ge 1}\sum_{\bk^{(1)}\!,\dots, \bk^{(s)}}\,\sum_{j_1+\dots+j_s=n-1}\prod_{r=1}^s\tfrac{t(j_r)}{j_r!}\prod_{u=1}^s\Bbb E[X_{\bk^{(u)}}(j_u)].
\end{eqnarray*}
Using $t(n)\Bbb E[X_{\bk}(n)]=a_{\bk}(n)$ with $a_{\emptyset}(n):=t(n)$ and $t(n)\Bbb E[Y_{\bk}(n)]=b_{\bk}(n)$, then multiplying the resulting equation by 
$\tfrac{z^{n-1}}{(n-1)!}$, and summing over $n\ge 2$, we  obtain
\begin{align*}
A_{\bk}'(z) & = \sum_{u\in [r]} B'_{k_u, (k_t)_{t\in [r]\setminus\{u\}}}(z)+\sum_{s\ge 1}\sum_{\bk^{(1)},\dots, \bk^{(s)}}\prod_{u=1}^s 
A_{\bk^{(u)}}(z),\\
&=\sum_{u\in [r]} B'_{k_u, (k_t)_{t\in [r]\setminus\{u\}}}(z)+\sum_{\tau\ge 1}(1-T(z))^{-\tau-1}\!\!\!\!\sum_{\bk^{(1)}\!,\dots, \bk^{(\tau)}\neq\emptyset}\,
\prod_{u=1}^{\tau} A_{\bk^{(u)}}(z).
\end{align*}
Explanation: The roots of $\tau$ subtrees that contain  non-empty sub-tuples of vertices of ranks $\bk^{(1)},\dots, \bk^{(\tau)}$ form $\tau+1$
consecutive intervals. Therefore, $(1-T(z))^{-\tau-1}$ is the EGF of finite forests of ordered trees whose roots are hosted by these intervals.

Let $\bk=(k_1,\dots,k_r)=:(k_1,\bk')$. To handle $B_{\bk}(z)$, introduce $B_{\ge k_1,\bk'}(z)$, the generating function of the numbers of occurrences of $\bv =(\text{root}, v_2,\dots, v_r)$\linebreak $=:(\text{root},\bv')$ such that the root is of rank $\ge k_1$, and $\bv'$ is of
rank $\bk'$.  Clearly $B_{\bold k}(z)\le B_{\ge k_1,\bk'}(z)$. Analogously to the case of $A_{\bk}(z)$, we obtain
\begin{equation} \label{new-17}
B_{\ge k_1,\bk'}'(z)\le \sum_{\tau\ge 1}\frac{1}{(1-T(z))^{\tau+1}}\biggl(\,\sum_{\cup_{u=1}^{\tau}\bk^{(u)}=\bk';\atop \bk^{(u)}\neq \emptyset}\,\,\prod_{u=1}^{\tau} A_{\bk^{(u)}}(z)\biggr).
\end{equation}
Here $\tau$ is the number of children of the root $1$ whose subtrees split among themselves the ordered $(r-1)$-tuple $\bk'$ according to a generic ordered partition of $\bk'$ into
$\tau$ non-empty sub-tuples $\bk^{(1)},\dots,\bk^{(\tau)}$, each ordered in compliance with the ordering of $\bk'$. Note that \eqref{new-17} is an inequality, because we
neglected the condition that the roots of all the subtrees, rooted at the children of the root $1$,  should have ranks $\ge k_1-1$ as well. 
%Furthermore, we have $\tau+1$ intervals that accommodate the remaining children of the root, subject to the condition the rank of each such child at its subtree is at least $k_1-1$.  Note that \eqref{new-17} is an inequality, because we
%neglected the condition that the roots of those $s$ trees should have ranks $\ge k_1-1$ as well. 
\end{proof}
}

%\subsection{Weighted sum of individual vertex counts by rank}\label{Sub1.2} Introduce $A(z,y)=\sum_{\ell\ge 0}y^{\ell}A_k(z)$, $B(z,y)=\sum_{\ell\ge 0}y^{\ell}B_k(z)$. Multiplying the equation for $A_k(z)$ by $y^k$ and adding we get
%\[
%A'_z(z,y)=\tfrac{A(z,y)}{(1-T)2}+B'_z(z,y).
%\]
 
\section{Asymptotics for the uniformly random tree $T_n$} \label{asymptotics}

We start this section by a series of lemmas that will allow us to prove an upper bound for the expected number
of vertices of rank at least $k$ in the uniformly random tree in Theorem \ref{thm1}. We then turn to the expected number of vertices of rank exactly $k$, and prove an asymptotic formula for that number in Theorem \ref{thm2}. 
That theorem will show that this number is asymptotically equal to $c_kn$, for a positive constant $c_k$.  We compute the values of $c_k$ for $k\leq 3$, and then, in Theorem \ref{add1}, we prove an upper bound for 
$1-\sum_{j\leq k} c_j$, which proves that the numbers $c_j$ form a probability distribution.

\subsection{Counting individual vertices by rank in $T_n$}
Recall that $A_k(z)$ is the (exponential) generating function of all occurrences of vertices of rank $k$ in all the trees on $[n]$,
$n\ge 1$, and $B_k(z)$ is the generating function of number of trees with root of rank $k$. We start with a simple, but
instrumental lemma.
\begin{lemma}\label{lem0}
For {\color{blue} $z\in \Bbb C$ and $|z|<1/2$}, we have
\[
A_k(z)=\frac{I_k(z)}{\sqrt{1-2z}},\quad I_k(z):=\int_0^z\sqrt{1-2\xi}\, B_k'(\xi)\,d\xi,
\]
{\color{blue} where the integral is taken over any path from $0$ to $z$ within the disc $|\xi|<1/2$.}
\end{lemma}
\begin{proof} According to Proposition \ref{general}, we have
\[
 A_k'(z)=\frac{A_k(z)}{(1-T(z))^2} +B_k'(z),\quad |z|<1/2.
 \]
 Since $T'(z)=(1-T(z))^{-1}$, we have $(1-T(z))^{-2}=T^{\prime\prime}(z)/T'(z)$, whence
 \[
 A'_k(z)=\frac{T^{\prime\prime}(z)}{T'(z)} A_k(z)+B'_k(z)\Longrightarrow \left(\frac{A_k(z)}{T'(z)}\right)'=\frac{B_k'(z)}{T'(z)}=
 (1-T(z))B_k'(z),
 \]
where $1-T(z)=\sqrt{1-2z}$. Our claim is now proved by integration.
 \end{proof}
%{
Let us introduce $A_{\ge k}(z)=\sum_{j\ge k} A_j(z)$, the exponential generating function for the number of all vertices of rank $\ge k$ in all rooted trees on vertex set $[n]$. By Lemma \ref{lem0}, we obtain
\begin{equation}\label{new-3}
A_{\ge k}(z)=\frac{I_{\ge k}(z)}{\sqrt{1-2z}},\quad I_{\ge k}(z):=\int_0^z\sqrt{1-2\xi}\, B_{\ge k}'(\xi)\,d\xi.
\end{equation}
%Now, the power series for $B_{\ge k}'(\xi)$ has only non-negative coefficients, and
%\[
%\sqrt{1-2\xi}=1-\sum_{\nu\ge 1}\frac{(2\nu-1)!!}{\nu!}.
%\]
%This implies that
%\[
%[\xi^{\mu}]

\begin{lemma}\label{lem1} For {\color{blue} $z\in \Bbb C$} and $|z|<1/2$ and  $k>2$, we have $|B'_{\ge k}(z)|\le 1/(k-2)!$. Consequently,
$|I_{\ge k}(z)|\le 1/(k-2)!$.
%\]
% \frac{(2z)^{k-2}}{(k-2)!}$. 
%Therefore
%$\lim_{k\to\infty}B'_{\ge k}(z)=0$ uniformly for $|z|<1/2$. Consequently
%\[
%c_{\ge k}:=\int_0^1\sqrt{1-2\xi}\, B_{\ge k}'(\xi)\,d\xi=O(1/(k-2)!).
%\]
\end{lemma} 
\begin{proof} 
Let {\color{blue} $z\in \Bbb R$} and $z\in [0,1/2)$. By Proposition \ref{general}, we have: for $\ell\ge 2$, 
%\begin{equation}\label{newer-1}
\begin{align*}
B_{\geq \ell}'(z)&=\frac{B_{\geq \ell-1}(z)}{1-B_{\geq \ell-1}(z)}\le\frac{B_{\ge \ell-1}(z)}{1-B_{\ge 1}(z)}\\
&=\frac{B_{\ge \ell-1}(z)}{1+z-T(z)}\le 2B_{\ge\ell-1}(z).
\end{align*}
{\color{blue} In the last step, we used the fact that if $z\in [0,1/2)$, then $z-T(z)\geq -1/2$.}
%\end{equation}
%Using the last identity and the main branch of the square root function, we analytically continue $T(z)$, whence $B_{\ge 1}(z)$, into an (open) domain $D_{\eps}:=\{z\in \Bbb C: |z|<0.5+\eps, z\notin [0.5,0.5+\eps)\}$. %Here $\rho=1+2^{3/2}/3$ is the root of $T(z)-z=1$ for $z>1/2$. 
%Keeping the notations $T(z)$ and $B_{\ge \ell}(z)$ for the continued versions, let us show that the above ODE allows
%to analytically continue each of $B_{\ge 2}(z), B_{\ge 3}(z),\dots $ into the domain $D_{\eps}$ for a sufficiently small
%$\eps>0$. 
%First of all, $\lim_{z\to 0.5}(T(z)-z) =1-1/2=0.5$, and $\sup_{z\in D_{\eps}}|T(z)-z|=0.5+O(\eps)$. Let $\ell\ge 2$ be such that, for each $j\in [1,\ell-1]$, $B_{\ge j}(z)$ is continued analytically to $D_{\eps}$ in such a way that $|B_{\ge j}(z)|\le 0.5 + O(\eps)$. This property holds for $\ell=2$. Then we continually extend $B_{\ge\ell}(z)$ from $|z|<1/2$ to $z\in D_{\eps}$
%by integrating the top ODE in \eqref{newer-2} and bound
%\begin{align*}
%|B_{\ge \ell}(z)|&\le \int_0^{|z|}\frac{|B_{\ge \ell-1}(y)|}{1-|B_{\ge \ell-1}(y)|}\,dy\le\frac{1}{0.5-O(\eps)}\int_0^{|z|}|B_{\ge \ell-1}(y)|\,dy\\
%&\le\cdots\le \frac{1}{(0.5-O(\eps))^{\ell-1}(\ell-2)!}\int_0^{|z|}(|z|-y)^{\ell-2}|B_{\ge 1}(y)|\,dy\\
%&\le \frac{(0.5+O(\eps))|z|^{\ell-1}}{(0.5-O(\eps))^{\ell-1}(\ell-1)!}\le \frac{(2+O(\eps))^{\ell-2}|z|^{\ell-1}}{(\ell-1)!}.
%\end{align*}

{\color{blue} Integrating} both sides of the resulting inequality $B_{\geq \ell}'(z)\le 2B_{\ge\ell-1}(z)$, by induction we obtain
\[
B_{\ge \ell}(z)\le \frac{2^{\ell-2}z^{\ell-1}}{(\ell-1)!}.
\]
Hence, for $k>2$, we have 
\[
B'_{\ge k}(z)\le 2B_{\ge k-1}(z)\le\frac{(2z)^{k-2}}{(k-2)!}\le \frac{1}{(k-2)!}.
\]
So, for {\color{blue} $z\in \Bbb C$ and $|z|<1/2$, we have}
\[
|B'_{\ge k}(z)|\le B'_{\ge k}(|z|)\le\frac{1}{(k-2)!},\quad |\sqrt{1-2z}|\le \sqrt{2}.
\]
It follows that $|I_{\ge k}(z)|\le 2^{-1/2} /(k-2)! \le 1/(k-2)!$.
\end{proof}

\noindent {\bf Notation.\/} To proceed, recall that $X_{k}(n)$ and $X_{\ge k}(n)$ stand for the total number of vertices of rank $k$
and of rank $\ge k$, respectively. Let $\mathcal R_n$ stand for the largest rank of a vertex in $T_n$.

\begin{theorem}\label{thm1}  Uniformly for  $k\le n$, we have
$\Bbb E[X_{\ge k}(n)]=O(n^{3/2}/(k-2)!)$. {\color{blue} Consequently, $\Bbb P\Bigl(\mathcal R_n\le (1.5+\eps)\frac{\log n}{\log\log n}\Bigr)=
1-O(n^{-\eps+o(1)})$.}
\end{theorem}
%{\bf Note.} Since the total number of vertices is $n$, it is the denominator $(k-2)!$ that makes the claim non-trivial. Whether in fact $\Bbb E[X_{n,\ge k}]=O(n/(k-2)!)$, or something similar, is an open problem.  

\begin{proof} By the Cauchy integral theorem, %and Lemma \ref{lem1}, 
for a circular contour $C$ of radius $\rho<1/2$ centered at zero, % $(C=\{z\in \Bbb C:\, z=\rho e^{i\theta},\, \theta\in [-\pi,\pi)\})$, 
we have
\[
\frac{a_{\ge k}(n)}{n!}=\frac{1}{2\pi i}\oint_{z\in C}\frac{I_{\ge k}(z)}{z^{n+1}\sqrt{1-2z}}\,dz.\\
%\le \bigl[\rho^n \sqrt{1-2\rho}\bigr]^{-1}\cdot\frac{1}{3(k-2)!}.
\]
{\color{blue} Let $C^*$ be the circle $|z|=1/2$.
 Then letting $C$ tend to $C^*$, substituting $z=\frac{1}{2} e^{i\theta}$, and using Lemma \ref{lem1}, we have
\begin{multline*}
\tfrac{a_{\ge k}(n)}{n!}\le \tfrac{2^{n+2}}{2\pi (k-2)!}\int_0^{\pi/2}|1-e^{i\theta}|^{-1/2}\,d\theta= \tfrac{2^{n+3/2}}{2\pi (k-2)!}\int_0^{\pi/2}
\sin^{-1/2}(\theta/2)\,d\theta\\
\le \tfrac{2^{n+2}}{2\pi (k-2)!}\int_0^{\pi/2}\theta^{-1/2}\,d\theta\le\tfrac{2^{n+1}}{(k-2)!}.
\end{multline*}
%$\frac{a_{\ge k}(n)}{n!}=O\bigl(2^n/(k-2)!\bigr)$, {\it uniformly for all $k\in [3,n]$\/}.
%The $\rho$-dependent factor attains its minimum $n^{1/2}(2+1/n)^{n+1/2}$ at $\rho=(2+1/n)^{-1}$, implying that
%\[
%\frac{a_{\ge k}(n)}{n!}\le\frac{(2e)^{1/2}}{(k-2)!}n^{1/2} 2^n.
%\]
So using $2^{2n}\binom{2n}{n}^{-1}\le 2n^{1/2}$, we get
%Therefore, using $m!=O(m^{1/2}(m/e)^m)$ and $(2n-1)!!=\Theta\bigl((2n/e)^n\bigr)$, we get
\begin{equation}\label{new-6.09}
\Bbb E[X_{\ge k}(n)]=\frac{a_{\ge k}(n)}{(2n-3)!!}\le \tfrac{2^{n+1}n!}{(2n-3)!! (k-2)! }\le \tfrac{4n}{(k-2)!}\cdot\frac{2^{2n}}{\binom{2n}{n}}\le 
\tfrac{8n^{3/2}}{(k-2)!}.
\end{equation}
Since $X_{k}(n)\le n$, this bound is useless if $k$ is fixed, or grows very slowly as $n$ grows. However,} setting $k_n=\lceil (1.5+\eps)\frac{\log n}{\log\log n}\rceil$, $\eps>0$ being arbitrary, and using $m!\ge (m/e)^m$, we conclude that
\[
\Bbb P(\mathcal R_n\ge k_n)\le \Bbb E[X_{\ge k_n}(n)]=O\bigl(n^{1.5} (k_n/e)^{-k_n}\bigr)=O\bigl(n^{-\eps'}\bigr),\quad \forall\,\eps'\in (0,\eps).
\]
\end{proof}
Let us turn now to $\Bbb E[X_{k}(n)]$. 
\begin{theorem}\label{thm2} For every fixed $k$, we have
\[
\Bbb E[X_{k}(n)]=c_k n \bigl(1+O(n^{-1})\bigr),\quad c_k:=2\int_0^{1/2}\!\!\! \!\sqrt{1-2\xi}\,B_k'(\xi)\,d\xi.
\]
\end{theorem}
\begin{proof}
By Lemma \ref{lem0}, for $|z|<1/2$,
\begin{equation} \label{New-5}
\begin{aligned}
&\quad\Bbb E[X_{n,k}]=\frac{a_k(n)}{(2n-3)!!} =\frac{n!}{(2n-3)!!} [z^n] A_k(z),\\
&A_k(z)=\frac{I_k(z)}{\sqrt{1-2z}},\quad I_k(z)=\int_0^z\!\!\sqrt{1-2\xi}\, B_k'(\xi)\,d\xi.
\end{aligned}
\end{equation}
Since $B_0(z)=z$, the integral $I_0(z)$ is $\bigl(1-(1-2z)^{3/2}\bigr)/3$, so we view $I_0(z)$ as being analytically extended to {\color{blue}$\Bbb C\setminus \{z\in \Bbb R: z\in [1/2,\infty)\}$}.

Further, recall that $B_k(z)=B_{\ge k}(z)-B_{\ge k+1}(z)$, and
\begin{equation}\label{new-6.1}
B_{\ge j+1}(z)=\int_0^z\frac{B_{\ge j}(\xi)}{1-B_{\ge j}(\xi)}\,d\xi,\quad |z|<1/2,
\end{equation}
{\color{blue} and the integral is taken along any path from $0$ to $z$ in the circle $|z|<1/2$. (Recall that $|B_j(z)|<1$ for $|z|<1/2$.) Simply by the definition of $B_{\ge 1}(z)$, we have
\[
B_{\ge 1}(z)=B_{\ge 0}(z) {\color{blue}-z }=\sum_{\nu\ge 2}\tfrac{t(\nu)}{\nu !} z^{\nu} =T(z)-z=1-\sqrt{1-2z}-z,
\]
with $\sqrt{\xi}=|\xi|^{1|/2} e^{i\text{Arg}(\xi)/2}$, $\text{Arg}(\xi)\in (-\pi,\pi)$. The last expression defines an analytic extension of $B_{\ge 1}(z)$ to
$\Bbb C\setminus [1/2,\infty)$.  For $|z|< 1/2$, the series formula implies that
$
|B_{\ge 1}(z)|\le B_{\ge 1}(|z|)\le B_{\ge 1}(1/2)=T(1/2)-1/2=1/2.
$
Since $B_{\ge j+1}(z)\le B_{\ge j}(z)$ for those $z$'s, we obtain that $|B_{\ge j}(z)|\le 1/2$ for $j>1$ and $|z|<1/2$.

So, by induction on $j$, it follows from \eqref{new-6.1} that for $|z|<1/2$ and  $j\ge 1$ the following holds.
\begin{enumerate} \item[(a)] $|B_j(z)|\le 1/2$, whence the integral in \eqref{new-6.1} is well-defined; 
\item [(b)]  $|B'_{\ge j+1}(z)|\le 1$.
\end{enumerate} 
%To be sure, if $|z|=1/2$, then the derivative is ``one-sided'', so that only $B_{\ge j+1}(z)$
%with $|z|\le 1/2$  are involved.
%if we define and uniformly for  points on the circle
%$|\xi|=1/2$, there exists  a finite $\ell_j=\lim_{z\to\xi: |z|<1/2}B'_{\ge j}(z)$ for each $j>1$.
We use  \eqref{new-6.1} to assert, by induction on $j$, the existence of  {\color{blue} the} analytic extension of each $B_{\ge j}(z)$ to the domain $D_j:=\{z\in \Bbb C: |z|<\rho_j\}\setminus [1/2,\rho_j)$, where $\rho_1\ge \rho_2\ge\cdots>1/2$, such that---
keeping $B_j(z)$ as notation for this extension---$|B_{\ge j}(z)|\le 3/4$ for $z\in D_j$, and \eqref{new-6.1} holds  for $z\in D_{j+1}$.
We use $D_j$ for {\it fixed\/} $j$ only. So, setting
$\rho>1/2$ equal to the minimum of $\rho_j$ over the {\it finitely\/} many $j$'s in question, we may and will focus on the analytic continuations $B_{\ge }j(z)$ to the same $D:=\{z\in \Bbb C: |z|<\rho\}\setminus [1/2,\rho)$.

Consequently, each $B_k(z)$ analytically extends to $D$.
 Keeping $B_k(z)$ and $I_k(z)$ to denote the extensions of $B_k(z)$
and the corresponding $I_k(z)$ given by \eqref{New-5}, we obtain {\color{blue} the following. 
Uniformly for }  $z\in D$, %and $k\ge 1$, 
we have
\begin{equation}\label{new-4.9}
\begin{aligned}
I_k(z)&=\int_0^z\sqrt{1-2\xi}\, B_k'(\xi)\,d\xi={\color{blue}\ga_k +\sum_{j\ge 0}\la_{k,j}(1-2z)^{\tfrac{j+2}{2}},}\\
\ga_k&:=\int_0^{1/2}\!\!\! \sqrt{1-2\xi}\,B_k'(\xi)\,d\xi;
\end{aligned}
\end{equation} 
 here the integral is taken over any path which is within $D$.  Recall that $B_1(z)=z$, so in particular, $B_1'(1/2)$ does exist. A bounded factor implicit in the big-O notation here, and  {\color{blue} in what follows}, may depend on $k$. 

Let $C_m$ be a contour that consists of the counter-clockwise circular arc $z=\rho e^{i\theta}$, $\theta\in (0,2\pi)$,
and a detour part formed by two opposite-directed line segments, one from $z=\rho e^{i((2-(1/m)\pi)}$ to $z=1/2 e^{i((2-(1/m))\pi)}$
and another from $z=1/2 e^{i(-1/m))}$ back to $z=\rho e^{i(1/m)}$.

Now consider the  contour $C$ that} consists of the counter-clockwise circular arc $z=\rho e^{i\theta}$, $\theta\in (0,2\pi)$,
and a detour part formed by two opposite-directed line segments, one from $z=\rho e^{i(2\pi-0)}$ to $z=1/2 e^{i(2\pi-0)}$
and another from $z=1/2 e^{i(+0)}$ back to $z=\rho e^{i(+0)}$. Then {\color{blue} $C$ is obviously a limit of the contours $C_m$, which are covered by the Cauchy theorem. 
See Figure \ref{contour} for an illustration.  (Readers can consult the Transfer theorem in 
\cite{flajolet} for more on this technique.)

\begin{figure}
 \begin{center}
  \includegraphics[width=120mm]{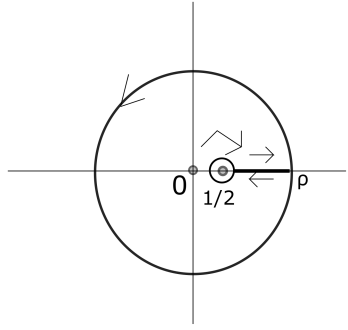}
  \caption{The contour $C$. }
  \label{fig:contour}
 \end{center}
\end{figure}

Applying this theorem and going to the limit, } we have, for $\a<1/2$,
\begin{multline*}%\label{new-4.91}
[z^n]\frac{I_k(z)}{\sqrt{1-2z}}=\frac{1}{2\pi i}\oint_{z\in \Bbb C: |z|=\a}\frac{I_k(z)}{z^{n+1}\sqrt{1-2z}}\,dz\\
=\ga_k [z^n] (1-2z)^{-1/2}+\la_{k,0} [z^n] (1-2z)^{1/2}+ \frac{1}{2\pi i}\oint_{z\in C}\frac{O\bigl(|1-2z|\bigr)}{z^{n+1}}\,dz.
\end{multline*}
Here
\begin{multline*}%\label{new-4.92}
 \ga_k[z^n] (1-2z)^{-1/2}+\la_{k,0} [z^n] (1-2z)^{1/2}=\ga_k\frac{(2n-1)!!}{n!}+\la_{k,0}\frac{(2n-3)!!}{n!}\\
=\ga_k\frac{(2n-1)!!}{n!}\bigl(1+O(n^{-1})\bigr)=\Theta\bigl(n^{-1/2}2^n\bigr).
\end{multline*}
and the remainder integral is of order
\begin{equation*}%\label{new-4.93}
\int_{1/2}^{\rho}\frac{2\eta-1}{\eta^n}\,d\eta +\frac{1}{\rho^n}\le n^{-2}2^n+O(\rho^{-n}).
\end{equation*}
We conclude that, for every fixed $k$,
\begin{equation}\label{new-4.94}
[z^n]\frac{I_k(z)}{\sqrt{1-2z}}=\ga_k\frac{(2n-1)!!}{n!}\bigl(1+O(n^{-1})\bigr).
\end{equation}
% So, by \eqref{New-5}, for every fixed $k$ we have
\[
\Bbb E[X_{k}(n)]=\frac{n!}{(2n-3)!!}\,[z^n]\frac{I_k(z)}{\sqrt{1-2z}}=(2n-1)\ga_k \bigl(1+O(n^{-1})\bigr),
\]
implying that $\Bbb E[X_{k}(n)]=nc_k\bigl(1+O(n^{-1})\bigr)$, $c_k:=2\ga_k$.
\end{proof}
%\begin{corollary}\label{corol2} By \eqref{new-6.09}, for each fixed $K$,
%\[
%n^{-1}\sum_{k\le K}\Bbb E[X_{n,k}]=1-n^{-1}\sum_{k>K}\Bbb E[X_{n,k}]\ge 1-\frac{5\sqrt{\pi e}}{(K-1)!}.
%\]
%Invoking Theorem \ref{thm2}, we have then that $1-\sum_{k\le K} c_k\le \frac{5\sqrt{\pi e}}{(K-1)!}$.
%\end{corollary}

\subsubsection{Computing $c_k$ for $k\leq 3$.} Recall that $c_k=\lim n^{-1}\Bbb E[X_{k}(n)]$, where
\[
c_k= 2\ga_k=2\int_0^{1/2}\!\!\! \sqrt{1-2\xi}\,B_k'(\xi)\,d\xi.
\]
Here $B_0(z)=z$, and for $k\ge 1$ we have the recurrence
\[
B_k(z)=B_{\ge k}(z)-B_{\ge k+1}(z),\quad B'_{\ge k}=\frac{B_{\ge k-1}(z)}{1-B_{\ge k-1}(z)},\quad |z|<1/2.
\] 
So,
\[
c_0=2\int_0^{1/2}\!\!\! \sqrt{1-2\xi}\,d\xi=\frac{2}{3}.
\]
Next, we compute
\[
B_1'(z)=B'_{\ge 1}(z)-B'_{\ge 2}(z)=\frac{B_{\ge 0}(z)}{1-B_{\ge 0}(z)}-\frac{B_{\ge 1}(z)}{1-B_{\ge 1}(z)};
\]
here $B_{\ge 0}(z)=T(z)=1-\sqrt{1-2z}$, and
\begin{equation}\label{new-6}
\begin{aligned}
B_{\ge 1}(z)=\int_0^z\frac{T(\xi)}{1-T(\xi)}\,d\xi&=\int_0^z\frac{1-\sqrt{1-2\xi}}{\sqrt{1-2\xi}}\,d\xi\\
&=1-z-\sqrt{1-2z}.
\end{aligned}
\end{equation}
Therefore
\begin{equation*}
B_1'(z)=\frac{1-\sqrt{1-2z}}{\sqrt{1-2z}}-\frac{1-z-\sqrt{1-2z}}{z+\sqrt{1-2z}}=\frac{z}{\sqrt{1-2z}\,(z+\sqrt{1-2z})}.
\end{equation*}
Consequently,
\begin{align*}
c_1&=2\int_0^{1/2}\!\!\! \sqrt{1-2z}\,B_1'(z)\,dz
=2\int_0^{1/2}\frac{z}{z+ \sqrt{1-2z}}\,dz\\%\biggl[\frac{1-\sqrt{1-2z}}{\sqrt{1-2z}}-\frac{1-z-\sqrt{1-2z}}{z+\sqrt{1-2z}}\biggr]\,dz\\
&=5+4\log2-6\sqrt{2}\,\text{arctanh} (2^{-1/2})\approx 0.2938858406.
\end{align*}
%with the numerical value delivered by Maple. %Further, integrating \eqref{new-6}, we have
%\begin{equation*}
%B_{\ge 1}(z)=-\frac{1}{3}+z-\frac{z^2}{2}+\frac{1}{3}(1-2z)^{3/2}.
%\end{equation*}
Next, 
\begin{eqnarray*}
B'_2(z) & = & B'_{\ge 2}(z)-B'_{\ge 3}(z)=\frac{1}{1-B_{\ge 1}(z)}-\frac{1}{1-B_{\ge 2}(z)}\\
& = & \frac{1}{z+\sqrt{1-2z}}-\frac{1}{1-B_{\ge 2}(z)},
\end{eqnarray*}
where 
%(Maple)
%\begin{align*}
%B_1(z)&=\int_0^z\frac{\xi}{\sqrt{1-2\xi}\,\bigl(\xi+\sqrt{1-2\xi}\bigr)}\,d\xi\\
%&=1-\sqrt{1-2z}-\log(z+\sqrt{1-2z})+2^{1/2}\text{arctanh}\,\bigl(2^{-1/2}(\sqrt{1-2z}-1)\bigr),
%\end{align*}
%and
%{\color{blue} 
\begin{eqnarray*}\label{new-7}
B_{\ge 2}(z) & = & \int_0^z\frac{B_{\ge 1}(y)}{1-B_{\ge 1}(y)}\,dy \\
& = & -z+\int_0^z\frac{1}{y+\sqrt{1-2y}}\,dy\\
%=-z+\log(z+\sqrt{1-2z})+2^{-1/2}\log\frac{1+2^{1/2}-\sqrt{1-2z}}{\sqrt{1-2z}-1+2^{1/2}}\\
& = & -z+\log(z+\sqrt{1-2z}) -2^{1/2}\text{arctanh}\bigl(2^{-1/2}(\sqrt{1-2z}-1)\bigr). 
\end{eqnarray*}

So, using the two last formulas, we obtain  {\color{blue}
\begin{eqnarray*}
c_2 & = & 2\int_0^{1/2}\sqrt{1-2z}\,B_2'(z)\,dz\\
& = & 2\int_0^{1/2}\sqrt{1-2z}\biggl(\frac{1}{z+\sqrt{1-2z}}-\frac{1}{1-B_{\ge 2}(z)}\biggr)\,dz\\
& \approx & 0.03589474655.
\end{eqnarray*} }
%\noindent $c_0+c_1+c_2\approx 0.66666666666+  0.2938858406+0.03589474655=0.9964472538$
Let us now turn to 
\begin{align*}
c_3&=2\int_0^{1/2} \sqrt{1-2z}\,B_3'(z)\,dz=2\int_0^{1/2} \sqrt{1-2z}\,\bigl(B_{\ge 3}'(z)-B_{\ge 4}'(z)\bigr)\,dz\\
&=2\int_0^{1/2} \sqrt{1-2z}\,\bigl((1-B_{\ge 2}(z))^{-1}-(1-B_{\ge 3}(z))^{-1}\bigr)\,dz.
\end{align*}
To compute $c_3$, introduce 
\[
C_3(z)=2\int_0^z \sqrt{1-2y}\,\bigl((1-B_{\ge 2}(y))^{-1}-(1-B_{\ge 3}(y))^{-1}\bigr)\,dy.
\]
Then $c_3=C_3(1/2)$ and $C_3(z)$, $B_{\ge 3}(z)$ satisfy a system of two ODEs,
\begin{align} \label{firstdeq}
\frac{dC_3}{dz}&=2\sqrt{1-2z}\,\bigl((1-B_{\ge 2}(z))^{-1}-(1-B_{\ge 3}(z))^{-1}\bigr),\\ \label{seconddeq}
\frac{dB_{\ge 3}}{dz}&=(1-B_{\ge 2}(z))^{-1}-1,
\end{align}
with $C_3(0)=B_{\ge 3}(0)=0$, and $B_{\ge 2}(z)$ being defined in the computation of $c_2$ above. 

%We used a recursive approximation based on the trapezoid method to compute $C_3(1/2)$ without computing $C_3(z)$ in general. 
%We obtain $c_3=C_3(1/2)= 0.0032684102$. This leads to 
%\[c_0+c_1+c_2+c_3\approx 0.9997156640.\]
 Sara Pollock, a  colleague of the first author, devised an efficient integration scheme of this ODE system, based on a trapezoidal approximation
of the integrals on each of intervals forming a partition of $[0,1/2]$, and used Matlab to show
that  $c_3:=C_3(1/2)=0.0032684102$ if each interval has length $h=10^{-6}$. Selecting $h=10^{-7}$ did not change the 
approximate value of $c_3$. We owe a debt of gratitude to Sara for her generous help.

%Thus $c_0+c_1+c_2+c_3\approx 0.9997156640$, and we already proved that $\sum_{k\ge 0} c_k=1$. So  the %limiting expected fraction of vertices having rank at least $4$ is $\approx 1- 0.9997156640=0.0002843360$.
{\color{blue}
\subsubsection{The numbers $\{c_k\}$ form a proper distribution.} 
As we stated in the introduction, we will prove that $\sum_{k\geq 0} c_k=1$, so the numbers $c_k$ form
a proper distribution. We saw in Section \ref{mainresults} that $\sum_{k\geq 0} c_k\leq 1$, so it suffices to show that
\[\lim_{k\rightarrow \infty}
\Bigl(1-\sum_{j\le k} c_j \Bigr)=0 .\] We accomplish this by proving the following stronger statement. 

\begin{theorem}\label{add1} For all $k$, the inequality
$
1-\sum_{j\le k} c_j\le \tfrac{3^{k+1}}{(2k+1)!}
$ holds.
\end{theorem} 
\noindent
{\bf Note.\/} So  the limiting expected fraction of vertices having rank at least $4$ is $0.0002843360\dots$ exactly.

\begin{proof} {\bf(i)\/}  Given  fixed $k$, let $a_{\le k}(n)$ and $a_{>k}(n)$ stand for the total number of 
{\color{blue} vertices of} rank $\le k$ and $>k$ in all trees on $[n]$, and let $X_{ \le k}(n)$, $X_{>k}(n)$ stand for the number of vertices of rank $\le k$ and $>k$ in the random tree.
So $t(n)\Bbb E[X_{\le k}(n)]=a_{\le k}(n)$ and $t(n)\Bbb E[X_{> k}(n)]=a_{>k}(n)$.  By \eqref{new-3}, 
\begin{equation}\label{nnew-3}
A_{\le k}(z)=(1-2z)^{-1/2}\int_0^z\sqrt{1-2x}\, B_{\le k}'(x)\,dx.
\end{equation}

The RHS of \eqref{nnew-3} is analytic within an open circle of a radius $\rho>1/2$ minus a cut $[1/2,\rho)$. It follows that 
\begin{align*}
\tfrac{(2n-3)!!\,\Bbb E[X_{\le k}(n)]}{n!} & = (1+o(1))[z^n](1-2z)^{-1/2}\int_0^{1/2}\sqrt{1-2x}\, B_{\le k}'(x)\,dx\\
& = (1+o(1))\tfrac{(2n-1)!!}{n!}\int_0^{1/2}\sqrt{1-2x}\, B_{\le k}'(x)\,dx\\
& = (1+o(1))\tfrac{(2n-1)!!}{n!}\int_0^{1/2}T'(x)B_{\le k}(x)\,dx.
\end{align*}

Using 
\begin{align*}
&\Bbb E[X_{\le k}(n)]=n-\Bbb E[X_{> k}(n)],\quad B_{\le k}(x)+B_{>k}(x)=T(x),\\
&c_j=\lim n^{-1}\Bbb E[X_{j}(n)], \quad 2\int_0^{1/2}T'(x) T(x)\,dx=1, 
\end{align*}
we obtain
\begin{equation}\label{Nnew-3}
1-\sum_{j\le k}c_j=2\int_0^{1/2}\tfrac{B_{>k}(x)}{1-T(x)}\,dx.
\end{equation}
Recall that,  for $|z|<1/2$, the power series  $B_{>k}(z)$ is the EGF generating function of $b_{>k}(n)$, where $b_{>k}(n)$ is the total number of trees in which the root is of rank $>k$. Clearly we need to find a good bound for $b_{>k}(n)$, or equivalently for $p_{>k}(n)=b_{>k}(n)/ t(n)$, the probability that the root of the random tree has rank more than $k$. To get such a bound we need a provably well-performing randomized algorithm which is likely to find a competitively short path from the root to a leaf. 

Here is the algorithm. If there are $s\ge 2$ children of the root, we keep a single tree, which we select with probability proportional to the total size of the deleted subtrees. For $n\ge 1$, $k\ge -1$, let $\pi_{>k}(n)$ denote the probability that edge length of the resulting path strictly exceeds $k$; obviously $p_{>k}(n)\le \pi_{>k}(n)$. Further, $\pi_{>- 1}(n)=1$, and for $n>1$, $k\ge 0$,
\begin{align*}
\pi_{>k}(n)& =\tfrac{t(n-1)}{t(n)}\,\pi_{>k-1}(n-1)\\
&\quad+  \tfrac{(n-1)!}{t(n)}\sum_{s\ge 2}\sum_{r=1}^s \sum_{j_1+\cdots+j_s=n-1}\!\!\tfrac{\sum_{r'\neq r}j_{r'}}{(n-1)(s-1)}\,\pi_{>k-1}(j_r)
\prod_{r"=1}^s\tfrac{t(j_{r"})}{j_{r"}!},
\end{align*}
or equivalently
\begin{align} \label{tnratio} 
 \notag \tfrac{t(n)\pi_{>k}(n)}{(n-2)!} &  =  \tfrac{t(n-1)\pi_{>k-1}(n)}{(n-2)!}  \\
& =  \sum_{s\ge 2}(s-1)^{-1}\sum_{r=1}^s\sum_{j_1+\cdots+j_s=n-1}\!\!\pi_{>k-1}(j_r)\sum_{r'\neq r}j_{r'}\prod_{r"=1}^s\tfrac{t(j_{r"})}{j_{r"}!}.
\end{align}
Introducing  $\Bbb P_{>k}(z)=\sum_{n\ge 1 }\tfrac{t(n)\pi_{>k}(n)}{n!}z^n$, (so that $\Bbb P_{>-1}(z)=T(z)$, and $\Bbb P_{>0}(z)=T(z)-z$),
we obtain a {\it linear\/} differential recurrence equation 
\begin{align*}\label{nnew-6}
\tfrac{d^2 \Bbb P_{>k}(z)}{dz^2} & =\tfrac{d \Bbb P_{>k-1}(z)}{dz}\\
+& \sum_{s\ge 2}\biggl(\sum_{j_1\ge 1}\tfrac{t(j_1)\pi_{>k-1}(j_1)}{j_1!}z^{j_1}\biggr)\biggl(\sum_{j_2\ge 1}\tfrac{t(j_2)}{(j_2-1)!} z^{j_2-1}\biggr)
\biggl(\sum_{j_3\ge 1}\tfrac{t(j_3)}{j_3!}z^{j_3}\biggr)^{s-2}  \notag     \\
=&  \tfrac{d \Bbb P_{>k-1}(z)}{dz}+\Bbb P_{>k-1}(z)\cdot T'(z)\sum_{t\ge 0}T^t(z), \notag
\end{align*}
or 
\begin{equation*}%\label{nnew-4}
\tfrac{d^2 \Bbb P_{>k}(z)}{dz^2}=\tfrac{d \Bbb P_{>k-1}(z)}{dz}+(1-T(z))^{-2}\Bbb P_{>k-1}(z),\quad |z|<1/2. 
\end{equation*}
Furthermore, by \eqref{Nnew-3}, for real-valued $x\in [0,1/2]$,
\begin{equation*}\label{nnew-5}
1-\sum_{j\le k}c_j\le 2\int_0^{1/2}\tfrac{\Bbb P_{>k}(x)}{1-T(x)}\,dx.
\end{equation*}
Since $(1-T(x))^{-1}=T'(x)$, we switch from $x$ to $\tau\in [0,1)$ such that $x=x(\tau)$ is inverse of $T(x)$. So, introducing  $\Pi_{>k}(\tau)=
\Bbb P_{>k}(x(\tau))$, we obtain
\begin{equation}\label{Nnew-5}
1-\sum_{j\le k}c_j\le 2\int_0^{1}\Pi_{>k}(\tau)\,d\tau.
\end{equation}
And, with a bit of work, the above differential equation becomes
\begin{equation}\label{nnew-4}
\Pi^{(2)}_{>k}(\tau)=\Pi^{(1)}_{>k-1}(\tau)\bigl(1-\tau-\tfrac{1}{1-\tau}\bigr)+\Pi_{>k-1}(\tau),\quad k\ge 0,\,\, \tau\in [0,1),
\end{equation}
subject to the initial conditions $\Pi_{>k}(0)=\Pi^{(1)}_{>k}(0)=0$.

Thus we need to bound the RHS of \eqref{nnew-5}. Introduce 
\begin{equation}\label{Ikt}
I_{k,t}=\int_0^{1}(1-\tau)^t\, \Pi_{>k}(\tau)\,d\tau, \quad k\ge -1,\,t\ge 0;
\end{equation}
so \eqref{Nnew-5} becomes
\[
1-\sum_{j\le k} c_j\le 2I_{k,0}.
\]
To bound $I_{k,0}$, we use \eqref{Ikt} to derive a recurrence inequality for all $t\ge 0$ and $k\ge -1$. First note that, for $t\ge 0$,
\begin{multline*}
I_{-1,t}=\int_0^1(1-\tau)^t\,\Pi_{>-1}(\tau)\,d\tau=\int_0^1 (1-\tau)^t \tau\,d\tau=\tfrac{1}{(t+2)_2}.\\
\end{multline*}
Let us turn to $k\ge 0$.
%Let us show that for $k\ge 0$ and $t\ge 0$
%\[
%I_{k,t}=
%\]
Integrating the RHS of \eqref{Ikt} by parts twice, using the initial conditions, and \eqref{nnew-4}, we have
\begin{eqnarray*}
I_{k,t} & = & \tfrac{1}{(t+2)_2}\int_0^1(1-\tau)^{t+2}\Pi_{>k}^{(2)}(\tau)\,d\tau\\
& = & \tfrac{1}{(t+2)_2}\biggl[\int_0^1\bigl[(1-\tau)^{t+3}-(1-\tau)^{t+1}\bigr]\Pi_{>k-1}^{(1)}(\tau)\,d\tau+I_{k-1,t+2}\biggr]\\
& = &  \tfrac{1}{(t+2)_2}\Bigl[(t+4)I_{k-1,t+2} - (t+1)I_{k-1,t}\Bigr]\\
& \le & \tfrac{3}{(t+2)_2}\cdot I_{k-1,t+2}.
\end{eqnarray*}
Here the third line results from yet another integration by parts, and the bottom line results   from the observation that $I_{k-1,s}$ decreases as $s$ increases.
It follows then that 
\[
I_{k,0}\le \tfrac{3^k}{(2k)!}\cdot I_{0,2k}\le \tfrac{3^k}{(2k+1)!}.
\]
\end{proof}

\subsubsection{How tight is $\{n^{-1}\Bbb E[X_{k}(n)]\}_{k\ge 0}$?} It trivially follows from Theorem \ref{add1} that the sequence of probability distributions $\{n^{-1}\Bbb E[X_{k}(n)]\}_{k\ge 0}$ is tight, meaning that $\sup_n n^{-1}\Bbb E[X_{>k}(n)]\to 0$ as 
$k\to\infty$. The question is how fast. The answer is super-exponentially fast.
\begin{theorem}\label{add2} For every $r \in (0,1)$, we have $\sup_n n^{-1}\Bbb E[X_{>k}(n)]=O(r^k)$, with the coefficient implicit  in the big-O notation dependent on $r$ only. %In words, the sequence $\{n^{-1}\Bbb E[X_{k}(n)]\}_{k\ge 0}$ is super-exponentially tight.
\end{theorem}
\begin{proof} {\bf (a)\/} Let $\rho>1$ and $\a\in (0,1)$. Introduce $g_n(\rho)=\sum_{k\ge 0}\rho^k \pi_{>k}(n)$. We want to prove that $g_n(\rho)\le h_n(\rho):=A n^{\a}$ if $A$ is sufficiently large. Let $n>1$. The, by \eqref{tnratio}, we have: 
\begin{eqnarray*}
g_n(\rho) & = & \tfrac{t(n-1)}{t(n)}(1+\rho g_{n-1}(\rho))\\
& + & \tfrac{(n-1)!}{t(n)}\sum_{s\ge 2}\tfrac{s}{(n-1)(s-1)} \sum_{j=1}^{n-s}(1+\rho g_j(\rho))\tfrac{t(j)(n-1-j)}{j!}\\
&\times& \sum_{j_1+\cdots+j_{s-1}=n-1-j}\prod_{r'=1}^{s-1}\tfrac{t(j_{r'})}{j_{r'}!}.
\end{eqnarray*}
Since $g_n(0)=\pi_{>0}(n)=1$ for $n>1$, and the $(s-1)$-fold sum equals $[z^{n-1-j}]\,T(z)^{s-1}$, the equation simplifies to 
\begin{align}\label{nnew-4.1}
g_n(\rho) & =  1+\tfrac{\rho t(n-1)}{t(n)}g_{n-1}(\rho)\\
& +  \tfrac{\rho\, (n-2)!}{t(n)}\sum_{s\ge 2}\tfrac{s}{s-1} \sum_{j=1}^{n-s}g_j(\rho)\tfrac{t(j)(n-1-j)}{j!}
\cdot [z^{n-1-j}] T(z)^{s-1}.
\end{align}
Changing the order of summation  in the double sum, we transform it into
\begin{eqnarray*}
\sum_{j=1}^{n-2}g_j(\rho)\tfrac{t(j)(n-1-j)}{j!}\cdot [z^{n-1-j}]\sum_{s=2}^{n-j}\tfrac{s}{s-1}T(z)^{s-1}\\
= \sum_{j=1}^{n-2}g_j(\rho)\tfrac{t(j)(n-1-j)}{j!}\cdot [z^{n-1-j}]\sum_{s=2}^{\infty}\tfrac{s}{s-1}T(z)^{s-1}\\
 = \sum_{j=1}^{n-2}g_j(\rho)\tfrac{t(j)(n-1-j)}{j!}\cdot [z^{n-1-j}]\biggl(\tfrac{T(z)}{1-T(z)}+\log\tfrac{1}{1-T(z)}\biggr)\\
=\sum_{j=1}^{n-2}g_j(\rho)\tfrac{t(j)}{j!}\cdot [z^{n-2-j}]\bigl(\tfrac{1}{(1-T(z))^3}+\tfrac{1}{(1-T(z))^2}\bigr);
\end{eqnarray*}
for the last step we used $\nu [z^{\nu}] F(z)=[z^{\nu-1}] F'(z)$ and $T'(z)=\tfrac{1}{1-T(z)}$. So \eqref{nnew-4.1} becomes
\begin{align} \label{nnew-4.2}
&g_n(\rho)=1+\tfrac{\rho t(n-1)}{t(n)}g_{n-1}(\rho)\\
&\qquad\quad+\tfrac{\rho\, (n-2)!}{t(n)}\sum_{j=1}^{n-2}g_j(\rho)\tfrac{t(j)}{j!}\cdot [z^{n-2-j}]\bigl(\tfrac{1}{(1-T(z))^3}+\tfrac{1}{(1-T(z))^2}\bigr).
%\times\sum_{s\ge 2}\tfrac{s}{s-1} \sum_{j=1}^{n-s}g_j(\rho)\tfrac{t(j)(n-1-j)}{j!}\cdot [z^{n-1-j}] T(z)^{s-1}.
\end{align}
Our task is to prove that the RHS of \eqref{nnew-4.2} where $g_j(\rho)$ is replaced with $h_j(\rho)$ is at most $h_n(\rho)$, {\it provided that
$n\ge n(\rho)$ and $A\ge A(\rho)$}. We will use
the identity
$
[z^{\nu}] (1-T(z))^{-\mu}=\tfrac{(\mu+2\nu-2)!!}{\nu!},
$
where the double factorial is the product of all positive integers at most $\mu+2\nu-2$, of the same parity as $\mu$.  Using the concavity of the function $x^{\a}$,
we have: $\sum_j p_j j^{\a}\le \bigl(\sum_j p_j j\bigr)^{\a}$ for $p_j\ge 0$, $\sum_jp_j=1$. 

To apply this inequality to
\[
\Sigma_n=\sum_{j=1}^{n-2}j^{\a}\,\tfrac{t(j)}{j!}\cdot [z^{n-2-j}]\bigl((1-T)^{-3}+(1-T)^{-2}\bigr),
\]
introduce 
\[
p_j:=\frac{\tfrac{t(j)}{j!}\cdot [z^{n-2-j}]\bigl((1-T)^{-3}+(1-T)^{-2}\bigr)}{\sum_{k}\tfrac{t(k)}{k!}\cdot [z^{n-2-k}]\bigl((1-T)^{-3}+(1-T)^{-2}\bigr)}
\]

The denominator equals
\[
[z^{n-2}]T(z)\bigl((1-T)^{-3}+(1-T)^{-2}\bigr)=[z^{n-2}]\, \bigl((1-T)^{-3}-(1-T)^{-1}\bigr).
\]

Therefore, for $n\ge 3$,
\begin{multline*}
%\sum_{j=1}^{n-2}j^{\a}\,\tfrac{t(j)}{j!}\cdot [z^{n-2-j}]\bigl((1-T)^{-3}+(1-T)^{-2}\bigr)\\
\Sigma_n\le \biggl([z^{n-2}]\, \bigl((1-T)^{-3}-(1-T)^{-1}\bigr)\biggr)^{1-\a}\\
\times\biggl(\sum_{j=1}^{n-2}\,\tfrac{jt(j)}{j!}\cdot [z^{n-2-j}]\bigl((1-T)^{-3}+(1-T)^{-2}\bigr)\biggr)^{\a},
\end{multline*}
Recalling $T'=(1-T)^{-1}$, we transforrn the second line sum into
\begin{multline*}
[z^{n-2}] \biggl(\sum_j\tfrac{jt(j)}{j!} z^j\biggr)\cdot \bigl((1-T)^{-3}+(1-T)^{-2}\bigr)\\
=[z^{n-2}]\,zT'\bigl((1-T)^{-3}+(1-T)^{-2}\bigr)=[z^{n-3}]\,T'\bigl((1-T)^{-3}+(1-T)^{-2}\bigr)\\
=[z^{n-3}]\,T'\bigl((1-T)^{-4}+(1-T)^{-3}\bigr).
\end{multline*}
Collecting the pieces, we obtain 
\begin{multline*}
\Sigma_n\le \biggl( [z^{n-2}]\, \bigl((1-T)^{-3}-(1-T)^{-1}\bigr)\biggr)^{1-\a}\\
\times\biggl([z^{n-3}]\bigl((1-T)^{-4}+(1-T)^{-3}\bigr)\biggr)^{\a}\\
=\Bigl(\tfrac{2(2n-5)!!}{(n-2)!}\Bigr)^{1-\a} \cdot \Bigl(\tfrac{(2n-4)!!+(2n-5)!!}{(n-3)!}\Bigr)^{\a}.
\end{multline*}
%\le [z^{n-2}]\, T(z)\bigl((1-T)^{-3}+(1-T)^{-2}\bigr)\\
%\times\Biggl(\frac{[z^{n-2}]\,zT'\bigl((1-T)^{-3}+(1-T)^{-2}\bigr)}{[z^{n-2}]\, T\bigl((1-T)^{-3}+(1-T)^{-2}\bigr)}\Biggr)^{\a}\\
%=\biggl( [z^{n-2}]\, \bigl((1-T)^{-3}-(1-T)^{-1}\bigr)\biggr)^{1-\a}\\
%\times\biggl([z^{n-3}]\bigl((1-T)^{-4}+(1-T)^{-3}\bigr)\biggr)^{\a}\\
%=\Bigl(\tfrac{2(2n-5)!!}{(n-2)!}\Bigr)^{1-\a} \cdot \Bigl(\tfrac{(2n-4)!!+(2n-5)!!}{(n-3)!}\Bigr)^{\a}.
%\end{multline*}
So the RHS of \eqref{nnew-4.2} with $h_j(\rho)$ instead of $g_j(\rho)$ is at most
\begin{multline*}
1+A\tfrac{\rho}{2n-3}(n-1)^{\a}+A\tfrac{\rho (n-2)!}{(2n-3)!!}\Bigl(\tfrac{2(2n-5)!!}{(n-2)!}\Bigr)^{1-\a} \cdot \Bigl(\tfrac{(2n-4)!!+(2n-5)!!}{(n-3)!}\Bigr)^{\a}\\
=1+A\tfrac{\rho}{2n-3}(n-1)^{\a}+A\rho\bigl(\tfrac{2}{2n-3}\bigr)^{1-\a}\Bigl((n-2)\bigl(\tfrac{(2n-4)!!}{(2n-3)!!}+\tfrac{1}{2n-3}\bigr)\Bigr)^{\a}\\
\le 1+ A\rho n^{-1+\a}+2A\rho n^{-1+2\a}\le 1+3A\rho n^{-1+2\a}.
\end{multline*}
The last bound is indeed below $An^{\a}$, if  $n>n(\rho)$ and $A\ge A(\rho,\a)$, where
\[
n(\rho):=1+\max\bigl\{3, (3\rho)^{(1-\a)^{-1}}\bigr\},\,\,A(\rho,\a):=\bigl(n^{\a}(\rho)-3\rho n^{2\a-1}(\rho)\bigr)^{-1}.
\]

Observe that, very crudely, $g_n(\rho)\le n_0 \rho^{n_0}$ if $n\le n_0$.  So, for those small values of $n$ we have  $g_n(\rho) < A n^{\a}$, if $A\ge n_0\rho^{n_0}$. By induction on $n\ge n_0$ we obtain  that, for all $n$, we have   $g_n(\rho)\le A^*n^{\a}$ if $A^*:=\max \{A(\rho,\a), n_0\rho^{n_0}\}$. Consequently, $p_{>k}(n) \le\pi_{>k}(n)\le A^*n^{\a}\rho^{-k}$.

\si
{\bf (b)\/}  %Introduce $H_n(\rho)=\sum_{k\ge 0}\rho^k \Bbb E[X_{>k}(n)]$. 
Using $h_n(\rho)\le An^{\a}$, we will prove that $\sup_{n\ge 1}n^{-1}\Bbb E[X_{>k}(n)]=O(r^k)$ for every $r<1$.  To 
begin, analogously to \eqref{??}, we have: for $n\ge 2$,
\begin{multline*}%\label{?}
\Bbb E[X_{>k}(n)]=p_{>k}(n)
+\tfrac{(n-1)!}{t(n)}\sum_{s\ge 1}\sum_{j_1+\dots+j_s=n-1}\prod_{r=1}^s\tfrac{t(j_r)}{j_r!}\biggl(\sum_{r'=1}^s
\Bbb E[X_{>k}(j_{r'})]\biggr).\\
\end{multline*}
By part {\bf (a)\/}, for all $\rho$ and $\a<1$, we have $\sum_k p_{>k}(n)\rho^k\le h_n(\rho)=An^{\a}$, with $A=A(\rho,\a)$.
Introduce $H_n(\rho)=\sum_{k\ge 0}\rho^k\Bbb E[X_{>k}(n)]$; obviously $H_n(\rho)=O(n)$ for $\rho\le 1$. Let us show
that $H_n(\rho)=O(n)$ for every $\rho>1$ as well. Of course, the coefficient hidden in the big-O notation will depend on $\rho$. From the above equation, and part {\bf (a)\/}, it follows that
\begin{multline*}
H_n(\rho)\le h_n(\rho)
+\tfrac{(n-1)!}{t(n)}\sum_{k\ge 0}\rho^k\sum_{s\ge 1}\sum_{j_1+\dots+j_s=n-1}\prod_{r=1}^s\tfrac{t(j_r)}{j_r!}\biggl(\sum_{r'=1}^s\Bbb E[X_{>k}(j_{r'})]\biggr)\\
=h_n(\rho)+\tfrac{(n-1)!}{t(n)}\sum_{j=1}^{n-1}\biggl(\sum_{k\ge 0}\rho^k\Bbb E[X_{>k}(j)]\biggr)\tfrac{t(j)}{j!}
\sum_{s\ge 1}s\sum_{\ell_1+\dots+\ell_{s-1}=n-1-j}\prod_{r=1}^{s-1}\tfrac{t(\ell_r)}{\ell_r!}\\
=h_n(\rho)+\tfrac{(n-1)!}{t(n)}\sum_{j=1}^{n-1}H_j(\rho)\tfrac{t(j)}{j!}\,[z^{n-1-j}] \sum_{s\ge 1}s T^{s-1}(z)\\
=h_n(\rho)+\tfrac{(n-1)!}{t(n)}\sum_{j=1}^{n-1}\,H_j(\rho)\tfrac{t(j)}{j!}[z^{n-1-j}] (1-T(z))^{-2}\\
=h_n(\rho)+\tfrac{2^{n-1}(n-1)!}{t(n)}\sum_{j=1}^{n-1}H_j(\rho)\tfrac{t(j)}{2^j j!},
\end{multline*}
since $(1-T(z))^{-2}=(1-2z)^{-1}=\sum_{u\ge 0}(2z)^u$. Clearly $H_n(\rho)\le \mathcal H_n(\rho)$, where $\mathcal H_1(\rho)=H_1(\rho)=0$ and
\[
\mathcal H_n(\rho)=h_n(\rho)+\tfrac{2^{n-1}(n-1)!}{t(n)}\sum_{j=1}^{n-1}\mathcal H_j(\rho)\tfrac{t(j)}{2^j j!},\quad n\ge 2.
\]
It follows that
\[
\tfrac{t(n+1)}{2^{n}n!}\bigl[\mathcal H_{n+1}(\rho)-h_{n+1}(\rho)\bigr]-\tfrac{t(n)}{2^{n-1}(n-1)!}\bigl[\mathcal H_{n}(\rho)-h_{n}(\rho)\bigr]=\mathcal H_n(\rho)\tfrac{t(n)}{2^n n!},
\]
or equivalently
\[ 
\tfrac{\mathcal H_{n+1}(\rho)}{2n+1}=\tfrac{\mathcal H_n}{2n-1} + D_n,\quad D_n:=\bigl[h_{n+1}(\rho)-\tfrac{2n}{2n-1}h_n(\rho)\bigr](2n+1)^{-1}.
\]
Since $h_n(\rho)=An^{\a}$, we have $D_n=O(n^{\a-1})$, and the recurrence equation implies that
\[
\mathcal H_{n+1}(\rho)=(2n+1)\cdot\biggl(\tfrac{\mathcal H_2(\rho)}{3}+\sum_{j=2}^n D_n\biggr)=O(n+n^{\a})=O(n).
\]
Therefore, uniformly for $k$, and $n\ge 1$,
\[
n^{-1}\Bbb E[X_{>k}(n)]\le n^{-1}\tfrac{\mathcal H_n(\rho)}{\rho^k}=O(\rho^{-k})
\]
for every $\rho>1$. In other words, the constant hidden in $O(\rho^{-k})$ depends only on $\rho$ and not on $k$.
\end{proof}

\subsection{Counting $r$-tuples of vertices with given ranks in $T_n$.} \label{rtuples}

%{\color{blue} 
We saw in  Proposition \ref{recur,gen}  that %with $\bk=\{k_1,\dots,k_r\}$
\begin{align}\label{new-20}
A'_{\bk}(z) & = \sum_{\tau\le r}(1-T)^{-\tau-1}\biggl(\,\sum_{\bk^{(1)}\!,\dots,\bk^{(\tau)}\neq \emptyset}\,\,\prod_{u=1}^{\tau} A_{\bk^{(u)}}(z)\biggr)\\
&\quad+\sum_{u\in [r]} B'_{k_u, (k_t)_{t\in [r]\setminus\{u\}}}(z), \notag
\end{align}
where $( \bk^{(1)},\dots,\bk^{(\tau)})$ denotes an ordered partition of an $r$-tuple $\bk$ into $\tau$ non-empty  sub-tuples, and
\begin{equation}\label{new-20.5}
B'_{\bk}(z)\le \sum_{\tau<r}\,\,(1-T)^{-\tau-1}\biggl(\,\sum_{\bk^{(1)}\!,\dots,\bk^{(\tau)}\neq \emptyset}\,\,\prod_{u=1}^{\tau} A_{\bk^{(u)}}(z) \biggr),
\end{equation}
where $\bk=(k_1,\bk')$ is an $r$-tuple, and $(\bk^{(1)},\dots, \bk^{(\tau)})$ denotes an ordered partition of $\bk'$ into $\tau$ non-empty sub-tuples.
%If $\bk'\neq\emptyset$, then
%\begin{align*}
%B'_{\bk}(z)&\le \sum_{\tau\ge 1}\,\,(1-T(z))^{-\tau-1}\biggl(\,\sum_{\bk^{(1)}\!,\dots,\bk^{(\tau)}\neq \emptyset}\,\,\prod_{u=1}^{\tau} A_{\bk^{(u)}}(z) \biggr);
%\end{aligned}
%\end{equation}
%in the equation $\{\bk^{(1)},\dots, \{bk^{(\tau)}\}$ denote a partition of $\bk$ into $\tau$ non-empty subtuples; in the inequality here $\bk=\{k_1,\bk'\}$, let $\{\bk^{(1)},\dots, \bk^{(\tau)}\}$ denote a generic ordered partition of $\bk'$ into $\tau$ non-empty ordered sub-tuples, and $\bk'\neq\emptyset$. 
Note that the right-hand side of \eqref{new-20.5}  does not depend on $k_1$.

%Thus the coefficients of the power series representing $A_{\bk}(z)$, i.e. the tree counts $a_{\bk}(n)$ divided by $n!$,
%are sandwiched between the coefficients for the power series solutions of the top equation in \eqref{new-20}, with $B'_{\bk}(z):\equiv 0$, and 
%\begin{equation}\label{new-20.1}
%B'_{\bk}(z):=\sum_{\tau\ge 1}\frac{1}{(1-T(z))^{\tau+1}}
%\sum_{\bk^{(1)}\,\dots, \bk^{(\tau)}\neq \emptyset}\,\prod_{u=1}^{\tau} 
%A_{\bk^{(u)}}(z).
%\end{equation}
%respectively.} To get an asymptotic formula for $[z^n] A_{\bk}(z)$ satisfying top line in \eqref{new-20}, it suffices to show that the contribution of $B'_{\bk}(z)$ defined by \eqref{new-20.1} to this coefficient is negligible.

\begin{lemma}\label{lem6} Given {\color{blue} an} $r$-tuple $\bold k$, there exist positive constants 
$\eps^{(i)}_{\bold k,\ell}$ such that  %there exists $\rho_{\bold k}>1/2$ such that for $|z|\le \rho_{\bk}$
we have $[z^{\nu}]A_{\bk}(z)\in \bigl[[z^{\nu}]A_{\bk}^{(1)}(z),[z^{\nu}]A_{\bk}^{(2)}(z)\bigr]$ for each $\nu$, ($A_{\bk}(z)\in [A_{\bk}^{(1)}(z),
A_{\bk}^{(2)}(z)]$, in short) where
\begin{equation}\label{new-21.4}
A_{\bk}^{(i)}(z)=\frac{(2r-3)!!\,\prod_{j=1}^r \ga_{k_j}}{(1-T)^{2r-1}}+\sum_{\ell\ge 2}\eps^{(i)}_{\bold k,\ell}(1-T)^{-2r+1+\ell},\quad T=T(z).
\end{equation}
\end{lemma}

}\begin{proof} This claim certainly holds for $r=1$; see \eqref{new-4.9}. Suppose the claim holds for all $r\le r_0$, for some $r_0\ge 1$. 
Each $\bk^{(u)}$ on the right-hand side  of the inequality \eqref{new-20.5} is of cardinality $\le r-1$. So, applying the inductive hypothesis to each $A_{\bk^{(u)}}(z)$ there, we bound the $B'$-sum in \eqref{new-20}:
\begin{equation}\label{new-21.5}
\sum_{u\in [r]} B'_{k_u, (k_t)_{t\in [r]\setminus\{u\}}}(z)\in \Bigl[ \sum_{\ell\le 2r-1}\delta^{(1)}_{\bold k,\ell}(1-T)^{-\ell}, \sum_{\ell\le 2r-1}\delta^{(2)}_{\bold k,\ell}(1-T)^{-\ell}\Bigr].
\end{equation}

{\color{blue}\noindent Let us turn to the $A$-sum in \eqref{new-20}. %$\prod_{u=1}^sA_{\bk^{(u)}}(z)$, $s\in [2,r]$. 
For $\tau\in [2,r]$, each $|\bk^{(u)}|\le r-1$, where $|\bk^{(u)}|$ is {\it just the number of components\/} in $\bk^{(u)}$.} So, by the inductive hypothesis again, the contribution of all ordered partitions of $[r]$ with $\tau\ge 2$ parts to this sum is at least
\begin{eqnarray*}
& & \frac{1}{(1-T)^{\tau+1}}\sum_{\cup_{u=1}^{\tau}\bk^{(u)}=\bk}
\,\,\prod_{u=1}^{\tau} A_{\bk^{(u)}}(z)\\
& = &   \frac{\prod_{j=1}^r\ga_{k_j}}{(1-T)^{2r+1}}\sum_{\cup_{u=1}^{\tau}\bk^{(u)}=\bk}
\prod_{u=1}^{\tau} \bigl(2|\bk^{(u)}|-3\bigr)!!\cdot\biggl(1+\sum_{\ell\ge 2}\eps_{\bk^{(u)},\ell}(1-T)^{\ell}\biggr)\\
& = &   \frac{\prod_{j=1}^r\ga_{k_j}}{(1-T)^{2r+1}}\sum_{\cup_{u=1}^{\tau}\bk^{(u)}=\bk}
\prod_{u=1}^{\tau} \bigl(2|\bk^{(u)}|-3\bigr)!! +\sum_{\ell\ge 2}\frac{\de^{(1)}_{\bk,s,\ell}}{(1-T)^{2r+1-\ell}}.
\end{eqnarray*}
And we have the similar bound from above.

%Here again each $|\bk^{(u)}|\le r-1$, since $s\ge 2$. By induction hypothesis, this contribution is 
%\[
%(1-T)^{-2r-1}\prod_{j=1}^r \ga_{k_j}+O(|1-T|^{-2r}),
%\]
%times the sum of the products $(2|\bk_j|-3)!!$ over all ordered partitions of $[r]$
%into $s$ parts.  
We need to sum these contributions over $\tau\in [2,r]$. The sums (those over $\ell\ge 2$) in the third line add up to the secondary sums
\[
\sum_{\ell\ge 2}\frac{\de^{(i)}_{\bk,\ell}}{(1-T)^{2r+1-\ell}},\quad \de^{(i)}_{\bk,\ell}:=\sum_{\tau=2}^r\de^{(i)}_{\bk,\tau,\ell},\quad i=1,2.
\]
The factor  $\frac{\prod_{j=1}^r\ga_{k_j}}{(1-T)^{2r+1}}$ aside, the first (primary) sum is the sum of the products 
$\prod_{u=1}^{\tau} \bigl(2|\bk^{(u)}|-3\bigr)!!$ over all ordered  partitions of the $r$ components of $\bk$, that is,  ordered partitions of the set $[r]$, into $\tau\ge 2$ non-empty subsets. We know that the total number of trees with $\be$ vertices is $(2\be-3)!!$. Therefore the sum of those products over
all ordered partitions of $\bk$ %over $s\in [1,r]$ 
is the total number of ordered
trees with $(r+1)$ vertices, (one root, some number of ordered trees, of total cardinality $r$, rooted at the root's children) i.e. $[2(r+1)-3]!!=(2r-1)!!$. The contribution coming from the partition with a single part of size $r$ is the total number of trees
with $r+1$ vertices (one root, one subtree with $r$ vertices rooted at a single child of the root), which is $(2r-3)!!$. Hence the sum of the products $\prod_{u=1}^{\tau}(2|\bk^{(u)}|-3)!!$ over all ordered partitions $(\bk_1,\dots,\bk_{\tau})$ of $[r]$
into $\tau\ge 2$ parts is $(2r-1)!! - (2r-3)!!=2(r-1)(2r-3)!!$. Consequently 
\[
\sum_{\cup_{u=1}^{\tau}\bk^{(u)}=\bk}
\prod_{u=1}^{\tau} \bigl(2|\bk^{(u)}|-3\bigr)!! =2(r-1)(2r-3)!!.
\]
So, according to \eqref{new-21.5}, it remains to solve the differential equation
\begin{align*}
\frac{dA^{(i)}_{\bk}(z)}{dz}=\frac{A^{(i)}_{\bk}(z)}{(1-T)^2}&+2(r-1)(2r-3)!! \cdot\frac{\prod_{j=1}^r \ga_{k_j}}{(1-T)^{2r+1}}\\
&+ \sum_{\ell\ge 2}\frac{\de^{(i)}_{\bk,\ell}}{(1-T)^{2r+1-\ell}}.
\end{align*}
%\[
%\frac{dA_{\bk}}{dz}=\frac{1}{(1-T)^2}A_{\bk}(z)+\sum_{s=2}^r\frac{s(2r-s-1)!}{2^{r-s}(r-s)!} \cdot\frac{\prod_{j=1}^r \ga_{k_j}}{
%(1-T)^{2r-1}}.
%\]
Using $T^{''}/T'=(1-T)^{-2}$, we obtain
\begin{eqnarray*}
\frac{d}{dz} \bigl[A_{\bk}(z)(1-T)\bigr]& =&2(r-1)(2r-3)!! \cdot\frac{\prod_{j=1}^r \ga_{k_j}}{(1-T)^{2r}}\\
&+& \sum_{\ell\ge 2}\frac{\de^{(i)}_{\bk,\ell}}{(1-T)^{2r-\ell}}.
\end{eqnarray*}
Recall that $1-T=(1-2z)^{1/2}$. So, upon integration from $0$ to $z$ in the open disc $|z|<1/2$, we have
\begin{equation*}
A_{\bk}(z)=\frac{(2r-3)!!\,\prod_{j=1}^r \ga_{k_j}}{(1-2z)^{r-1/2}}+\sum_{\ell\ge 2}\eps_{\bold k,\ell}(1-2z)^{-r+1/2+\ell/2},
\end{equation*}
which is the equation \eqref{new-21.4}. 
\end{proof}
\begin{theorem}\label{thm3} Given an $r$-tuple $\bk=(k_1,\dots, k_r)$, let $X_{\bold k}(n)$ denote the total number of $r$-tuples of
distinct vertices $(v_1,\dots, v_r)$ in the random tree $T_n$ with ranks $\bk$. Then
\[
\Bbb E[X_{\bk}(n)]=(1+O(n^{-1})) n^r\prod_{j=1}^r c_{k_j}.
\]
\end{theorem}
\begin{proof} By the definition of $A_{\bk}(z)$, we have
\[
\Bbb E[X_{\bk}(n)]=\frac{n! [z^n] A_{\bk}(z)}{(2n-3)!!}.
\]
To evaluate this coefficient, introduce $C_{\bk}$, a contour enclosing $z=0$, which is analogous to $C$ in \eqref{new-4.9}-\eqref{new-4.94}, with $\rho:=\rho_{\bk}$. By Lemma \ref{lem6}, and following the derivation in those equations, we obtain
\begin{align*}
 [z^n] A_{\bk}(z) 
=  (2r-3)!!\,\biggl(\prod_{j=1}^r\ga_{k_j}\biggr)[z^n](1-2z)^{-r+1/2}\\
 +  \sum_{\ell=2}^{2r}\eps_{\bk,\ell}\, [z^n] (1-2z)^{-r+1/2+\ell/2}+
\frac{1}{2\pi i}\oint_{z\in C_{\bk}}\frac{O\bigl(|1-2z|\bigr)}{z^{n+1}}\,dz\\
=  (2r-3)!!\biggl(\prod_{j=1}^r\ga_{k_j}\biggr)\frac{\bigl(2(n+r)-3\bigr)!!}{n!(2r-3)!!}\\
+  O\biggl(\frac{(2(n+r)-5)!!}{n!}\sum_{\ell=2}^{2r}\eps_{\bk,\ell}\biggr) +O\biggl(\int_{1/2}^{\rho_{\bk}}\frac{(2u-1)}{u^{n+1}}\,du +\rho_{\bk}^{-n}\biggr)\\
=\bigl(1+O(n^{-1})\bigr)\biggl(\prod_{j=1}^r\ga_{k_j}\biggr)\frac{\bigl(2(n+r)-3\bigr)!!}{n!}\biggr)+O(n^{-2}2^n).
\end{align*}
We can drop the last error term, because the explicit term is of order $n^{r-3/2} 2^n$ exactly. So, since $c_k=2\ga_k$. 
\begin{align*}
\Bbb E[X_{\bk}(n)]&=\frac{n!\cdot [z^n] A_{\bk}(z)}{(2n-3)!!}=(1+O(n^{-1}))\biggl(\prod_{j=1}^r\ga_{k_j}\biggr)\cdot
\prod_{u=1}^r \bigl(2(n+u)-3\bigr)\\
&=(1+O(n^{-1}))\,n^r\prod_{j=1}^r c_{k_j}.
\end{align*}
\end{proof}
\begin{corollary}\label{corol3} Let $R_n(1),\dots, R_n(r)$ denote the ranks of $r$ vertices chosen uniformly at random, with order and without replacement, from $[n]$. Then, for every  $(k_1,\dots, k_r)$, we have
\[
\Bbb P(R_n(1)=k_1,\dots, R_n(r)=k_r)=(1+O(n^{-1})) \prod_{j=1}^r c_{k_j}.
\]
In words, the ranks of a finite, uniformly random, set of vertices are asymptotically independent, each with distribution
$\{c_k\}_{k\ge 0}$.
\end{corollary}
%{\color{blue}{\bf Note.\/} We haven't proved yet that the distribution $\{c_k\}_{k\ge 0}$ is {\it proper\/}, i.e. $\sum_{k\ge 0}c_k=1$. So, for now, the limits of $R(j)$ assume infinite values with deficit probability $1-\sum_{j\ge 0}c_j$.}

%
\subsection{More on the largest rank $\mathcal R_n$.} {\color{blue} In Section 3.1 we proved that with probability $\ge 1-n^{- 0.99\eps}$,  $\mathcal R_n$ is
at most $(1.5+\eps)\log n/\log\log n$. Here we will prove that  with probability $\ge 1-n^{- 0.99\eps}$,  the inequality $\mathcal R_n\ge (1-\eps)\log n/\log\log n$ also holds.} 

We focus on vertices $v$ for which the  subtree rooted at $v$ is a path connecting $v$ to a leaf. For example, vertices whose only child is a leaf are vertices with this property. For brevity, we call such vertices path-type, or p-type. The rank of a p-type vertex $v$ is the edge-length of the path that connects $v$ to a unique leaf. Therefore $\mathcal R_n\ge \mathcal R^p_n$, where $\mathcal R_n^p$ is the largest rank of a path-type vertex in $[n]$.

Recall that $A_k(z)$ is the exponential generating function for the number of all vertices of rank $k$ in all rooted
plane trees on vertex set $[n]$, and $B_k(z)$ is the exponential generating function for the number of trees in which the root is of rank $k$. Furthermore, recall $A_{(k,k)}(z)$ is the exponential generating function for the number of 
ordered pairs  of vertices of rank $k$ each  in all trees on $[n]$, and $B_{(k,k)}(z)$ is the exponential generating function for the number of ordered pairs $(v,w)$ of rank $k$ each in all such trees so that $v$ is the {\em root} of the tree.
We will need the following lemma. 

\begin{lemma}\label{lem7}  Let $\mathcal A_k(z)$, $\mathcal A_{(k,k)}(z)$, $\mathcal B_k(z)$, and $\mathcal B_{(k,k)}(z)$
be the counterparts of  $A_k(z)$, $A_{(k,k)}(z)$, $B_k(z)$, and $B_{(k,k)}(z)$ for p-type vertices. Then
$\mathcal B_k(z)=z^k/k!$, while $\mathcal B_{(k,k)}(z)=0$, and  consequently,
\begin{equation}\label{new-21.51}
\begin{aligned}
 \mathcal A_{k}'(z)=&\frac{\mathcal A_{k}(z)}{(1-T(z))^2}+\frac{z^{k-1}}{(k-1)!},%\mathcal B'_{k}(z),
 \\
%\mathcal B_{k}'(z) =& \mathcal B_{k-1}(z), \,\, k\ge 1;\,\,\, \mathcal B_{0}(z)=z,
%\end{aligned}
%\end{equation}
%and
%\begin{equation}\label{new-21.52}
%\begin{aligned}
\mathcal A_{(k,k)}'(z)=&\frac{\mathcal A_{(k,k)}(z)}{(1-T(z))^2}+\frac{2\mathcal A_{k}^2(z)}{(1-T(z))^3}.
%+
%\mathcal B'_{\bk}(z),
%\mathcal B_{\bk}'(z)=&
\end{aligned}
\end{equation}
\end{lemma} 
%\begin{proof} 
\noindent The proof is a simpler version of the proofs of Proposition \ref{recur} and Proposition \ref{recur,gen},  and we omit it. 
Note that the part concerning $\mathcal B_k(z)$ and $\mathcal B_{(k,k)}(z)$ is straightforward. Indeed, there is only one increasing tree with a p-type root of rank $k$, namely the path $1\to 2\to\cdots\to k$. So $\mathcal B_k(z)$, the exponential generating function of numbers of trees with a p-type root of rank $k$,
is given by $\mathcal B_k(z)=z^k/k!$. And there are no trees such that both the root and another vertex are
p-type and are of the same rank $k$. So $\mathcal B_{(k,k)}(z)$, the exponential generating function of numbers of trees with the root and another vertex both p-type and of the same rank $k$, is zero.

Given $k$, let $\mathcal X_{k}(n)$ denote the total number of path-type vertices in the random tree $T_n$ of rank $k$.
\begin{theorem}\label{thm4}  Given $\eps\in (0,1)$, define $k=[(1-\eps)(\log n)/\log\log n]$. Then
\[
\frac{\Bbb E[\mathcal X^2_{ k}(n)]}{\Bbb E^2[\mathcal X_{k}(n)]}=1+O(\eps).
\]
Consequently,  with probability $\ge 1- O(n^{-\eps})$,  there exist $\Theta(n^{\eps})$ vertices (in fact, p-type vertices) with rank $k(n)$.
\end{theorem}

Before we can prove Theorem \ref{thm4}, we need to 
solve the system of  ODEs that appear in Lemma \ref{lem7}. Arguing as in the proof of Lemma \ref{lem0}, we have
\begin{equation}\label{new-21.52}
\begin{aligned}
\mathcal A_k(z)&=\frac{1}{\sqrt{1-2z}}\int_0^z\sqrt{1-2\xi}\,\,\frac{\xi^{k-1}}{(k-1)!}\,d\xi,\\
\mathcal A_{(k,k)}(z)&=\frac{1}{\sqrt{1-2z}}\int_0^z\frac{2\mathcal A_k^2(\xi)}{1-2\xi}\,d\xi.
\end{aligned}
\end{equation}
This leads to the following Lemma that we will use to prove Theorem \ref{thm4}. 

\begin{lemma}\label{lem8} For $k<n$, we have
\begin{equation}\label{new-21.53}
\begin{aligned} \
[z^n] \mathcal A_k(z)&=\frac{(2n-1)!!}{n!\,(2k+1)!!},\\
[z^n] \mathcal A_{(k,k)}(z)&=\bigl(1+O(n^{-1})\bigr)\,\frac{(2n+1)!!}{n!\,[(2k+1)!!]^2}.
%\mathcal A_k(z)&=\sum_{j\le k}\frac{(-1)^j\binom{k}{j}}{2^k\,k!\,(2j+3)}\Bigl[(1-2z)^{-1/2}-(1-2z)^{j+1}\Bigr],\\
\end{aligned}
\end{equation}
\end{lemma}
\begin{proof} Let us start with $\mathcal A_k(z)$. We compute
\begin{align*}
\sqrt{1-2\xi}\,\,\frac{\xi^{k-1}}{(k-1)!}&=\frac{\sqrt{1-2\xi}}{2^{k-1}\,(k-1)!}\bigl(1-(1-2\xi)\bigr)^{k-1}\\
&=\frac{1}{2^{k-1}\,(k-1)!}\sum_{j<k}(-1)^j\binom{k-1}{j}(1-2\xi)^{j+1/2}.
\end{align*}
Therefore, integrating
\begin{align*}
\mathcal A_k(z)&=\frac{(1-2z)^{-1/2}}{2^{k-1}\,(k-1)!}\sum_{j<k}\frac{(-1)^j}{2j+3}\binom{k-1}{j}\bigl[1-(1-2z)^{j+3/2}\bigr]\\
&=\frac{(1-2z)^{-1/2}}{2^{k-1}\,(k-1)!}\sum_{j<k}\frac{(-1)^j}{2j+3}\binom{k-1}{j}-P_k(z), \\
P_k(z)&:=\frac{1}{2^{k-1}\,(k-1)!}\sum_{j<k}\frac{(-1)^j}{2j+3}\binom{k-1}{j} (1-2z)^{j+1}.
\end{align*}
Here $P_k(z)$ is a polynomial of degree $k$. Let us assume from now on that $k<n$.  Then
\[
[z^n]\,\mathcal A_k(z)=\frac{1}{2^{k-1}\,(k-1)!}\sum_{j<k}\frac{(-1)^j}{2j+3}\binom{k-1}{j}\times [z^n](1-2z)^{-1/2}.
\]
The familiar last factor is $(2n-1)!!/n!$, and (using Maple in the third line) \footnote{Note that 
$\int_0^1 x^2 (1-x^2)^{k-1} \ dx $ equals, by the change of variables $u=x^2$, the expression $  B(3/2,k)/2$, where $B$ is the Beta-function. Then using the formula for Beta in terms of the Gamma-function the formula follows.}

\begin{eqnarray*}
\sum_{j<k}\frac{(-1)^j}{2j+3}\binom{k-1}{j} & = & \sum_{j<k}(-1)^j\binom{k-1}{j}\int_0^1 x^{2j+2}\,dx\\
& = & \int_0^1\!x^2(1-x^2)^{k-1}\,dx \\
& = & \int_0^1\!(1-x^2)^{k-1}\,dx- \int_0^1(1-x^2)^k\,dx\\
& = & \frac{2^{k-1}\,(k-1)!}{(2k+1)!!}. 
\end{eqnarray*}
We conclude that 
\begin{equation}\label{new-21.54}
%\begin{aligned}
\mathcal A_k(z)=\frac{(1-2z)^{-1/2}}{(2k+1)!!}-P_k(z),\quad [z^n]\mathcal A_k(z)=\frac{(2n-1)!!}{n!\,(2k+1)!!}.
%\end{aligned}
\end{equation}
It remains to plug the formula for $\mathcal A_k(z)$ into the bottom equation in 
\eqref{new-21.52}. It is easy to check that neglecting $P_k(z)$ in this formula results in an asymptotic formula
for $[z^n] \mathcal A_{(k,k)}(z)$ within a factor $1+O(n^{-1})$ from the actual coefficient. With $P_k(z)$ dropped, we have
\begin{align*}
[z^n] (1-2z)^{-1/2}\!\int_0^z\frac{2(1-2\xi)^{-2}}{[(2k+1)!!]^2}\,d\xi 
& =\!  \frac{1}{[(2k+1)!!]^2} [z^n]\Bigl[\tfrac{1}{(1-2z)^{-/2}}-\tfrac{1}{(1-2z)^{1/2}}\Bigr]\\
&= \!\frac{1}{n! [(2k+1)!!]^2}\Bigl[(2n+1)!! -(2n-1)!!\Bigr] \\
& =\!  (1+O(n^{-1}))\frac{(2n+1)!!}{n!\,[(2k+1)!!]^2}.
\end{align*}
Therefore
\begin{equation}\label{new-21.55}
[z^n] \mathcal A_{(k,k)}(z)=\bigl(1+O(n^{-1})\bigr)\,\frac{(2n+1)!!}{n!\,[(2k+1)!!]^2}.
\end{equation}
\end{proof}

Now we can prove Theorem \ref{thm4}. 
\begin{proof} (of Theorem \ref{thm4}). 
Using \eqref{new-21.54}, %and \eqref{new-21.55} 
we have
\[
\Bbb E[\mathcal X_{k}(n)]=\frac{n! [z^n]\mathcal A_k(z)}{(2n-3)!!}=\frac{(2n-1)}{(2k+1)!!}.
\]
Applying Stirling formula,  and setting $k(n):=\bigl[\tfrac{(1-\eps)\log n}{\log\log n}\bigr]$, we obtain
\[
\Bbb E[\mathcal X_{k(n)}(n)]=n^{\eps}\bigl(1+O(1/(\log\log n))\bigr)\to\infty.
\]
Furthermore
\[
\Bbb E[(\mathcal X_{k(n)}(n))_2]=\frac{n! [z^n]\mathcal A_{k(n),k(n)}(z)}{(2n-3)!!}
=(1+O(n^{-1}))\Bbb E^2[\mathcal X_{k(n)}(n)].
\]
So,
\begin{align*}
\Bbb E[\mathcal X^2_{k(n)}(n)]&=\Bbb E[\mathcal X_{k(n)}(n)] +\Bbb E[(\mathcal X_{k(n)}(n))_2]\\
&=\bigl(1+O(n^{-1}) +O(\Bbb E^{-1}[\mathcal X_{k(n)}(n)])\bigr)\Bbb E^2[\mathcal X_{k(n)}(n)]\\
&=(1+O(n^{-\eps}))\Bbb E^2[\mathcal X_{k(n)}(n)].
\end{align*}
Consequently, by Chebyshev's inequality, we have
\[
\Bbb P\Bigl[\Big|\frac{\mathcal X_{k(n)}(n)}{\Bbb E[\mathcal X_{k(n)}(n)]}-1\Big|\ge \delta]=
O(\delta^{-2}n^{-\eps}).
%\Bbb P\Bigl(\mathcal X_{n,k}\ge 0.5 \Bbb E[\mathcal X_{n,k}]\Bigr)\ge 1- O(n^{-\eps}).
\]
\end{proof}

Combining Theorem \ref{thm1} and Theorem \ref{thm4}, we arrive at the following result. 
\begin{theorem}\label{thm5} For the largest rank $\mathcal R_n$ of a vertex, and $\eps\in (0,1)$,
\[
\Bbb P\biggl(\frac{\mathcal R_n}{\frac{\log n}{\log\log n}}\in [1-\eps, 1.5+\eps]\biggr)=1-O\bigl(n^{-0.99\eps}\bigr).
\]
\end{theorem}
\noindent {\color{blue}Thus the largest vertex rank slowly but inexorably tends to infinity. Still the ranks of finitely many, uniformly random vertices approach finite random limits since $\{c_j\}$ is proper.  Can $1.5$ be replaced with $1$ without changing the RHS bound?}
%random limits since
%\begin{theorem}\label{thm6} $\sum_{j\ge 0} c_j=1$.
%\end{theorem}
%\begin{proof} Let $V_k(n)$ denote the total number of vertices whose descendant trees are of size $k$. Let $\mathcal C_k(z)$ denote the
%EGF of the sequence $\{t(n)\Bbb E[V_k(n)]\}$; then, analogously to \eqref{new-21.51},
%\[
%\mathcal C'_k(z)=\frac{\mathcal C_k(z)}{(1-T)^2}+\bigl(\tfrac{t(k)}{k!} z^k\bigr)'.
%\]
%so that $\mathcal C_k(z)= t(k)\mathcal A_k(z)$. Using \eqref{new-21.53}, we have
%\begin{align*}
%\tfrac{t(n)\Bbb E[V_k(n)]}{n!}&=t(k)[z^n]\mathcal A_k(z)=\tfrac{(2n-1)!!\, t(k)}{n! (2k+1)!!}\\
%&\quad\Longrightarrow \Bbb E[V_k(n)]=\tfrac{(2k-3)!!\,(2n-1)!!}{(2k+1)!!\,(2n-3)!!}=\tfrac{2n-1}{(2k-1)(2k+1)}.
%\end{align*}
%Therefore
%\[
%d_k:=\lim_{n\to\infty}\tfrac{\Bbb E[V_k(n)]}{n}=\tfrac{2}{(2k-1)(2k+1)}=\tfrac{1}{2k-1}-\tfrac{1}{2k+1},
%\]
%and $\sum_{k\ge 1}d_k=1$. Thus $\{d_k\}$ is a proper probability distribution of the descendant tree size which is rooted at a uniformly random vertex. Now, for every $K$, $\sum_{k\le K+1}d_k$ is the limiting probability that the tree size is at most $K+1$, and on this event the rank of tree's root is at most $K$. Hence 
%\begin{equation}\label{room}
%\sum_{k\le K}c_k\ge \sum_{k\le K+1}d_k=1-\tfrac{1}{2K+3}.
%\end{equation}
%Letting $K\to\infty$, we obtain $\sum_{k\ge 0}c_k=1$.
%\%end{proof}
%{\bf Note.\/} Our numerical results indicate that it should be possible to drastically reduce the remainder term $O(K^{-1})$ to something exponentially small, say.
\bi
{\bf Acknowledgment.\/} We are genuinely grateful to the two hard-working referees who helped us to significantly improve our paper.

\bi
{\color{blue} {\bf Appendix.\/} Originally we used the randomized algorithm for a suboptimal path from the root to a leaf for the search tree, \cite{bona-pittel}.
This is a plane binary tree on vertex set $[n]$: recursively, vertex $n+1$ joins the left subtree with probability proportional to its current size. We
proved the existence of the limiting distribution $\{c_j\}$ of $R_n$, the rank of the random vertex in $[n]$, satisfying 
\[
1-\sum_{j=0}^{k-1}c_j\le \tfrac{6k+7}{3}\bigl(\tfrac{1}{3}\bigr)^k.
\]
We also claimed that for every $\rho\in (0, 3/2)$, the equality $\lim \Bbb E[\rho^{R_n}]=\sum_j \rho^j c_j$ holds, implying that $R_n$ converges, in distribution, to the limit rank $\mathcal R_{\infty}$ with all its moments. Looking at the proof anew, we have found an oversight. Here we present a corrected argument that proves a stronger result:
$\lim \Bbb E[\rho^{R_n}]=\sum_j \rho^j c_j$ for every $\rho\in [0,3)$, so that 
\[
1-\sum_{j=0}^{k-1}\Bbb P(R_n=j) =O\bigl((\tfrac{1}{3}+\eps)^k\bigr),
\]
for every $\eps>0$, almost matching the above tail bound for $\{c_j\}$. The proof is based on the randomized algorithm. Let $L_n$ and  $\mathcal L_n$ denote, respectively, the rank of the root and the edge length of the path from the root to a leaf delivered by the algorithm. 
Then the inequality $\Bbb P(L_n>k)\le \Bbb P(\mathcal L_n>k)$ holds, meaning that $h_n(\rho):=\Bbb E[\rho^{L_n}]\le \Bbb E[\rho^{\mathcal L_n}]$.
From the definition of  plane binary trees and that of the algorithm,  $\pi_{n,>k}:=\Bbb P(\mathcal L_n>k)$ satisfies
$\pi_{n,>-1}=1$ for $n\ge 1$, and for $k\ge 0$ we have $\pi_{1,>k}=0$, while for  $n>1$ $\pi_{n,>k}$ satisfies a recurrence equation
\begin{align*}
\pi_{n,>k}&=\tfrac{1}{n}\sum_{j=0}^{n-1} \Bigl(\tfrac{n-1-j}{n-1}\, \pi_{j,>k-1}+\tfrac{j}{n-1}\,\pi_{n-1-j,>k-1}\bigr)\\
&= \tfrac{2}{n(n-1)}\sum_{j=0}^{n-1}(n-1-j)\pi_{j,>k-1};\quad \pi_{0,>-1}:=1.
%2(n-1)\pi_{n-1,>k-1}+2\sum_{j=1}^{n-2}(n-1-j)\pi_{j,>k-1},\quad n\ge 1,\,\,k\ge -1.
\end{align*}
%\begin{equation*}
%n(n-1)\pi_{n,>k}=2(n-1)\pi_{n-1,>k-1}+2\sum_{j=1}^{n-2}(n-1-j)\pi_{j,>k-1},\quad n\ge 1,\,\,k\ge -1.
%\end{equation*}
So, for $n>1$,  $g_n(\rho):=\sum_{k\ge -1}\rho^k\pi_{n,>k}$ satisfies
\begin{equation}\label{a-1}
g_n(\rho)=\rho^{-1}+\sum_{k\ge 0} \rho^k \pi_{n,>k}
=\rho^{-1}+\tfrac{2}{n(n-1)}\sum_{j=0}^{n-1}(n-1-j)g_j(\rho).
\end{equation}
Let us select $\rho>1$ and $\a\in (0,1)$. Since $g_j(\rho)<j\rho^j$, for every $J>0$ there exists $A:=A(\a,\rho,J)>0$ such that $g_j(\rho)\le Aj^{\a}$ for $j\le J$. Now, 
for $n>J$ we have
\[
\sum_{j=1}^{n-2}(n-1-j)j^{\a}\le n\int_0^ny^{\a}\,d\a -\int_0^ny^{\a+1}\,dy=\tfrac{n^{\a+2}}{(\a+1)(\a+2)}.
\]
So, the RHS of \eqref{a-1} with $Aj^{\a}$ instead of $g_j(\rho)$ is at most
\[
\rho^{-1}+A\tfrac{2\rho n^{\a}}{(\a+1)(\a+2)} (1+O(J^{-1})),
\]
which in turn is below $An^{\a}$ for $J$ large enough, provided that $\rho<\tfrac{(\a+1)(\a+2)}{2}$ and
\[
A>\frac{\rho^{-1}}{1-\tfrac{2\rho}{(\a+1)(\a+2)}(1+O(J^{-1}))}.
\]
Notice that $\rho$ can be pushed  arbitrarily close from below to $\rho=3$ by making $\a$ close enough to $1$ from below. By induction on
$n\ge J$ we obtain that $g_n(\rho)\le A n^{\a}$ for all $n$.

Consequently, $h_n(\rho)=\Bbb E[\rho^{L_n}]= O(n^{\a})$ for $\rho<\tfrac{(\a+1)(\a+2)}{2}$. Introduce $H_n(\rho)=\sum_{k<n}\rho^k\Bbb E[X_{n,k}]$, the expected value of $\sum_{v\in [n]}\rho^{R(v)}$, where $R(v)$ denotes the rank of a vertex $v$. Then
\[
H_n(\rho)=h_n(\rho)+\frac{1}{n}\sum_{j=0}^{n-1}\bigl(H_j(\rho)+H_{n-1-j}(\rho)\bigr),\quad n>1.
\]
Here $\tfrac{h_n(\rho)}{n}=O(n^{\a-1})$; applying Lemma 2.1 from \cite{bona-pittel} we obtain: for each $\rho<\tfrac{(\a+1)(\a+2)}{2}$, there exists a finite 
$
\lim_{n\to\infty} n^{-1}\sum_{v\in [n]}\rho^k \Bbb E[X_{n,k}].
$
In other words, denoting $R_n$ the rank of the random vertex in the random search tree, we have: for each $\rho<3$,
$
\lim_{n\to\infty} \sum_{k\ge o} \rho^k\, \Bbb P(R_n=k)=\sum_{k\ge 0}\rho^k c_k.
$
%\bi
%{\bf Note.\/} It seems plausible that a similar claim holds for the plane tree in the present paper as well.
}
\bi

\end{document}